\documentclass[12pt]{amsart}
\usepackage{stmaryrd}
\usepackage{pifont}
\usepackage{amsmath, mathtools}
\usepackage{mathrsfs}
\usepackage{bbm}
\usepackage{amsfonts}
\usepackage{latexsym,amssymb,amsmath,amscd,amsfonts,epsfig,graphicx}

\usepackage[margin=1.4in]{geometry}

\newtheorem{thm}{Theorem}[section]
\newtheorem{lem}[thm]{Lemma}

\newtheorem{coro}[thm]{Corollary}
\newtheorem{prop}[thm]{Property}%[section]
\newtheorem{propo}[thm]{Proposition}%[section]
\newtheorem{obs}[thm]{Observation}%[section]

\def\pf{\noindent{\it Proof.\;\;\:}}
\def\qed{\nopagebreak\hfill{\rule{4pt}{7pt}}
\medbreak}

\DeclareMathOperator{\diam}{diam}
\DeclareMathOperator{\radi}{radi}
\DeclareMathOperator{\dist}{dist}
\DeclareMathOperator{\ecc}{ecc}

\DeclareMathOperator{\girth}{girth}
\theoremstyle{definition}
\newtheorem{defi}[thm]{Definition}%[section]
\newtheorem{exa}[thm]{Example}%[section]

\title[Original graphs of link graphs]{Original graphs of link graphs}

\author[Bin Jia]{Bin Jia}

\thanks{Bin Jia gratefully acknowledges scholarships provided by The University of Melbourne.}

\begin{document}

\maketitle

\begin{center}
Department of Mathematics and Statistics\\
The University of Melbourne\\
Victoria 3010, Australia.\\
Email: jiabinqq@gmail.com \quad Mobile: +61 404639816
\end{center}

\medskip

\noindent {\bf Abstract.}
Let $\ell \geqslant 0$ be an integer, and $G$ be a graph without loops.
An \emph{$\ell$-link} of $G$ is a walk of length $\ell$ in which consecutive edges are different. We identify an $\ell$-link with its reverse sequence. The \emph{$\ell$-link graph $\mathbb{L}_\ell(G)$} of $G$ is defined to have vertices the $\ell$-links of $G$, such that two vertices of $\mathbb{L}_\ell(G)$ are adjacent if their corresponding $\ell$-links are the initial and final subsequences of an $(\ell + 1)$-link of $G$. A graph $G$ is called an \emph{$\ell$-root} of a graph $H$ if $\mathbb{L}_\ell(G) \cong H$. For example, $\mathbb{L}_0(G) \cong G$. And the $1$-link graph of a simple graph is the line graph of that graph. Moreover, let $H$ be a finite connected simple graph. Whitney's isomorphism theorem (1932) states if $H$ has two connected nonnull simple $1$-roots, then $H \cong K_3$, and the two $1$-roots are isomorphic to $K_3$ and $K_{1, 3}$ respectively.

This paper investigates the $\ell$-roots of finite graphs. We show that every $\ell$-root is a certain combination of a finite minimal $\ell$-root and trees of bounded diameter. This transfers the study of $\ell$-roots into that of finite minimal $\ell$-roots. As a qualitative generalisation of Whitney's isomorphism theorem, we bound from above the number, size, order and maximum degree of minimal $\ell$-roots of a finite graph. This work forms the basis for solving the recognition and determination problems for $\ell$-link graphs in our future papers. As a byproduct, we characterise the $\ell$-roots of some special graphs including cycles.

Similar results are obtained for path graphs introduced by Broersma and Hoede (1989). $G$ is an \emph{$\ell$-path root} of a graph $H$ if $H$ is isomorphic to the $\ell$-path graph of $G$. We bound from above the number, size and order of minimal $\ell$-path roots of a finite graph.

%Tree-decompositions and better-quasi-orderings of the original graphs of finite link and path graphs are also investigated. As a byproduct, we generalise a theorem of Ding (1992) as follows: for any infinite sequence of finite or infinite graphs (respectively, of bounded multiplicity), if none of them contains an $\ell$-path as a subgraph, then there are two graphs in the sequence such that the former is a (respectively, an induced) subgraph of the latter. Moreover, in this case, these graphs admit linked tree-decompositions of minimal tree-width such that the trees are of bounded diameter.

\medskip


\section{Introduction and main results}
As a generalisation of line graphs \cite{Whitney1932} and path graphs \cite{BH1989}, the link graph of a graph $G$ was introduced by Jia and Wood \cite{JiaWood2013} who studied the connectedness, chromatic number and minors of the link graph based on the structure of $G$. This paper deals with the reverse; that is, for an integer $\ell \geqslant 0$ and a finite graph $H$, we study the graphs whose $\ell$-link graphs are isomorphic to $H$.

Unless stated otherwise, all graphs are undirected and loopless. A graph may be finite or infinite, and may be simple or contain parallel edges. In particular, $H$ always denotes a finite graph. The \emph{order} and \emph{size} of $H$ are \emph{$n(H) := |V(H)|$} and \emph{$m(H) := |E(H)|$} respectively. Throughout this paper, $\ell \geqslant 0$ is an integer. An \emph{$\ell$-link} is a walk of length $\ell$ in which consecutive edges are different. We identify an $\ell$-link with its reverse sequence. An \emph{$\ell$-path} is an $\ell$-link without repeated vertices. The \emph{$\ell$-link graph $\mathbb{L}_\ell(G)$} of a graph $G$ is defined to have vertices the $\ell$-links of $G$, and two vertices are adjacent if their corresponding $\ell$-links are the initial and final subsequences of an $(\ell + 1)$-link of $G$. Furthermore, if $G$ contains parallel edges, then by \cite[Observation 3.4]{JiaWood2013}, two $\ell$-links of $G$ may be the initial and final subsequences of each of $\mu \geqslant 2$ different $(\ell + 1)$-links of $G$. In this case, there are $\mu$ edges between the two corresponding vertices in $\mathbb{L}_\ell(G)$. Rigorous definitions are given in Section \ref{sec_term}.

A graph $G$ is called an \emph{$\ell$-root} of $H$ if $\mathbb{L}_\ell(G) \cong H$. Let $\mathbb{R}_\ell(H)$ be the set of minimal (up to the subgraph relation) $\ell$-roots of $H$. The \emph{line graph $\mathbb{L}(G)$} of $G$ is the simple graph with vertex set $E(G)$, in which two vertices are adjacent if their corresponding edges share a common end vertex in $G$. By definition $\mathbb{L}_1(G) = \mathbb{L}(G)$ if and only if $G$ is simple. So Whitney's isomorphism theorem \cite{Whitney1932} can be restated as this: Let $H$ be a connected nonnull simple graph. Then $|\mathbb{R}_1(H)| \leqslant 2$, with equality holds if and only $H \cong K_3$ and $\mathbb{R}_1(H) = \{K_3, K_{1,3}\}$. Moreover, it is not difficult to see that every $1$-root of $K_3$ is the disjoint union of one of $K_3$ and $K_{1, 3}$, and zero or more isolated vertices. Motivated by Whitney's isomorphism theorem, this paper aims to answer the following three general questions for a finite graph $H$ and an integer $\ell$ to some extent:
\begin{itemize}
\item How many minimal $\ell$-roots can $H$ have?
\item How large can a minimal $\ell$-root of $H$ be?
\item How to construct all $\ell$-roots of $H$ from minimal $\ell$-roots of $H$?
\end{itemize}

We answer the first two questions by the following theorem, which is a qualitative generalisation of Whitney's isomorphism theorem.

\begin{thm}\label{thm_minirootbounded} Let $\ell \geqslant 0$ be an integer, and $H$ be a finite graph. Then
the maximum degree, order, size, and total number of minimal $\ell$-roots of $H$ are finite and bounded by functions of $\ell$ and $|V(H)|$.
\end{thm}
\noindent {\bf Remark.} The order and size of minimal $\ell$-roots of $H$ are bounded in Lemma~\ref{lem_HboundG}. Then other parameters are trivially bounded. But we are able to improve the remaining bounds by further investigating the structure of link graphs. In particular, $|\mathbb{R}_\ell(H)|$ is bounded in Lemma \ref{lem_lminiFinite}. And the maximum degree of a minimal $\ell$-root is bounded in Corollary \ref{coro_DsdecreaseCyc} and Lemma \ref{lem_Tdc}.

Before continuing, we give more motivations for the first question. First of all, for some infinite families of graphs $H$, $|\mathbb{R}_\ell(H)|$ is bounded by a constant number. For example, Lemma \ref{lem_cyclerootT} characterises all minimal $\ell$-roots of a cycle, which implies that the number of minimal $\ell$-roots of a cycle is one or two. However, $|\mathbb{R}_\ell(H)|$ is not always bounded by a constant number. For instance, Lemma \ref{lem_barK2} finds all minimal $\ell$-roots of $2K_1$, which indicates that the number of minimal $\ell$-roots of $2K_1$ increases with $\ell$. On the other hand, in Section \ref{sec_example} we exemplify that, for fixed $\ell \geqslant 4$ and any given number $k$, there exists a connected graph $H$ with $k$ or more minimal connected $\ell$-roots. The last two examples show that in Theorem \ref{thm_minirootbounded}, it is necessary that the number of minimal $\ell$-roots of $H$ are bounded by functions of $\ell$ and $|V(H)|$.

%
%Tree-decompositions of $\ell$-roots are worth considering. Let $G$ be a hypergraph, $T$ be a tree, and $\mathcal{V} := \{V_w | w \in V(T)\}$ be a set cover of $V(G)$ indexed by the nodes of $T$. The pair $(T, \mathcal{V})$ is called a \emph{tree-decomposition} of $G$ if for each edge $e$ of $G$, there exists some $V \in \mathcal{V}$ containing all vertices incident to $e$, and for every path $[w_0, \ldots, w_i, \ldots, w_\ell]$ of $T$, we have $V_{w_0} \cap V_{w_\ell} \subseteq V_{w_i}$. The \emph{width} of the tree-decomposition $(T, \mathcal{V})$ is defined to be \emph{$\tw(T, \mathcal{V}) := \sup\{|V| - 1| V \in \mathcal{V}\}$}. The \emph{tree-width $\tw(G)$} is the minimum width of a tree-decomposition of $G$. The \emph{tree-diameter $\tdi(G)$} is the minimum diameter of $T$ over the tree-decompositions $(T, \mathcal{V})$ of $G$ of width $\tw(G)$. The following statement is proved in Section \ref{sec_pathGraphs}.
%
%\begin{thm}\label{thm_llinkbqo} Let $\ell \geqslant 0$ be an integer, and $H$ be a finite graph. Then
%the tree-width and tree-diameter of the $\ell$-roots of $H$ are finite and bounded by functions of $|V(H)|$ and $\ell$. Furthermore, the $\ell$-roots of a finite graph are better-quasi-ordered by the induced subgraph relation.
%\end{thm}

Denote by \emph{$X \subseteq Y$}, \emph{$X \subset Y$}, \emph{$X \leqslant Y$}, and \emph{$X < Y$} that $X$ is isomorphic to a subgraph, proper subgraph, induced subgraph and proper induced subgraph of a graph $Y$ respectively. A graph is \emph{$\ell$-finite} if its $\ell$-link graph is finite. So all finite graphs are $\ell$-finite, but not vice versa. For example, let $T^{t}$ be the tree obtained by pasting the middle vertex of a $4$-path at the center of a star $K_{1, t}$. Then $T^{\infty}$ is infinite, and is $4$-finite since its $4$-link graph is $K_1$.

Two $\ell$-finite graphs
$X$ and $Y$ are \emph{$\ell$-equivalent}, write \emph{$X \thicksim_\ell Y$}, if there exists a graph $Z \subseteq X, Y$ such that $\mathbb{L}_\ell(X) \cong \mathbb{L}_\ell(Y) \cong \mathbb{L}_\ell(Z)$. For example, for every pair of integers $i, j \geqslant 0$, we have $T^i \thicksim_4 T^j$ since $\mathbb{L}_4(T^i) \cong \mathbb{L}_4(T^j) \cong \mathbb{L}_4(T^0) \cong K_1$.
An $\ell$-finite graph $X$ is \emph{$\ell$-minimal} if $X$ is null or $\mathbb{L}_\ell(Y) \subset \mathbb{L}_\ell(X)$ for every $Y \subset X$. For instance, an $\ell$-path is $\ell$-minimal. By definitions, a graph is $\ell$-minimal if and only if it is a minimal $\ell$-root of a finite graph.

Let \emph{$\mathbb{R}_\ell[H]$} be the set of $\ell$-roots of $H$. By the analysis above, $\mathbb{R}_4[K_1]$ is an equivalence class with equivalence relation $\thicksim_4$. Moreover, $\mathbb{R}_1[K_3]$ is the union of two equivalence classes with equivalence relation $\thicksim_1$. And the two classes contain $K_3$ and $K_{1, 3}$ respectively. The lemma below is proved in Section \ref{sec_consEquiClass}. Together with Theorem \ref{thm_minirootbounded}, it says that $\thicksim_\ell$ divides $\mathbb{R}_\ell[H]$ into finitely many equivalence classes.

\begin{lem}\label{lem_miniUniq}
For each integer $\ell \geqslant 0$, $\sim_\ell$ is an equivalence relation on $\ell$-finite graphs, such that each $\ell$-equivalence class contains a unique (up to isomorphism) $\ell$-minimal graph. And this graph is isomorphic to an induced subgraph of every graph in its class.
\end{lem}

The \emph{distance $\dist_G(u, v)$} between $u, v \in V(G)$ is $+ \infty$ if $u, v$ are in different components of $G$, and the minimum length of a path of $G$ between $u, v$ otherwise. The \emph{eccentricity} of $v \in V(G)$ is \emph{$\ecc_G(v) := \sup\{\dist_G(u, v) | u \in V(G)\}$}.  The \emph{diameter} of $G$ is \emph{$\diam(G) := \sup\{\ecc_G(v)| v \in V(G)\}$}.

Lemma \ref{lem_conClass} answers the third general question of this paper, which is also proved in Section \ref{sec_consEquiClass}. It says that an $\ell$-root of a finite graph is a certain combination of a minimal $\ell$-root and trees of bounded diameter. This transfers the study of $\ell$-roots into that of minimal $\ell$-roots. In Section \ref{sec_example}, we first obtain all minimal $\ell$-roots of a cycle. Then we characterise all $\ell$-roots of a cycle by applying the lemma below.

\begin{lem}\label{lem_conClass}
Let $G$ be the minimal graph of an $\ell$-equivalence class. Then a graph belongs to this class if and only if it can be obtained from $G$ in two steps:
\begin{itemize}
\item[{\bf (1)}] For each acyclic component $T$ of $G$ of diameter within $[\ell, 2\ell - 4]$, and every vertex $u$ of eccentricity $s$  in $T$ such that $\lceil\ell/2\rceil \leqslant s \leqslant \ell - 2$, paste to $u$ the root of a rooted tree of height at most $\ell - s - 1$.
\item[{\bf (2)}] Add to $G$ zero or more acyclic components of diameter at most $\ell - 1$.
\end{itemize}
\end{lem}

%The theorem below says that the order, size, maximum degree, and total number of minimal $\ell$-roots are finite and bounded.

%Lemma \ref{lem_conClass} will be applied in Section \ref{sec_example} to obtain all $\ell$-roots of a cycle.
% transfers the study of $\ell$-roots into that of minimal $\ell$-roots whose

%Lemma \ref{lem_conClass} transfers the study of $\ell$-roots into that of minimal $\ell$-roots. We have mentioned that the order, size, maximum degree, and total number of the $\ell$-roots of a given graph are unbounded. In contrast, our main result below confirms that all these parameters are bounded for minimal $\ell$-roots.

%{thm_minirootbounded}

Introduced by Broersma and Hoede \cite{BH1989}, the \emph{$\ell$-path graph $\mathbb{P}_\ell(G)$} is the simple graph with vertices the $\ell$-paths of $G$, where two vertices are adjacent if the union of their corresponding paths forms a path of length $\ell + 1$ or a cycle of length $\ell + 1$ in $G$. Let $\ell \geqslant 2$. It follows from the definition that $\mathbb{P}_\ell(G)$ is a subgraph of $\mathbb{L}_\ell(G)$. And the two graphs are isomorphic if and only if $\girth(G) > \ell$. For each $\ell \geqslant 0$, we say $G$ is an \emph{$\ell$-path root} of $H$ if $\mathbb{P}_\ell(G) \cong H$. Let $\mathbb{Q}_\ell(H)$ be the set of minimal (up to the subgraph relation) $\ell$-path roots of $H$. Li \cite{LiXL1996} proved that $H$ has at most one simple $2$-path root of minimum degree at least $3$. Prisner \cite{Prisner2000} showed that $\mathbb{Q}_\ell(H)$ contains at most one simple graph of minimum degree greater than $\ell$.
By Li and Liu \cite{LiLiu2008}, if $H$ is connected and nonnull, then $\mathbb{Q}_2(H)$ contains at most two simple graphs. In fact, the finite graphs having exactly two minimal simple $2$-path roots have been characterised by Aldred, Ellingham, Hemminger and Jipsen \cite{AEHJ1997}. Some results about $\ell$-roots can be proved, with slight variations, for $\ell$-path roots:

\begin{thm}\label{thm_pathrootsBound} Let $\ell \geqslant 0$ be an integer, and $H$ be a finite graph. Then
the order, size, and total number of minimal $\ell$-path roots of $H$ are finite and bounded by functions of $\ell$ and $|V(H)|$.
\end{thm}

\section{Terminology}\label{sec_term}
This section presents some definitions and simple facts. The reader is referred to \cite{Diestel2010} for notation and terminology on finite and infinite graphs. A graph is said to be \emph{cyclic} if it contains a cycle, and \emph{acyclic} otherwise. Let $G$ be a graph, and \emph{$c(G)$} (respectively, \emph{$o(G)$}, \emph{$a(G)$}) be the cardinality of the set of (respectively, cyclic, acyclic) connected components of $G$. A \emph{ray} is an infinite graph with vertex set $\{v_0, v_1, \ldots\}$ and edges $e_i$ between $v_{i - 1}$ and $v_i$, for $i \geqslant 1$.
The \emph{radius $\radi(G)$} of $G$ is $+\infty$ if $G$ is disconnected, and $\min\{\ecc_G(v)| v \in V(G)\}$ otherwise. For each tree $T$, Wu and Chao \cite{WuChao2004} proved that $\radi(T) = \lceil \diam(T)/2\rceil$.

Denote by \emph{$K_t$} the complete graph on $t \geqslant 0$ vertices. In particular, \emph{$K_0$} is the \emph{null} graph. Let $tG$ be the disjoint union of $t \geqslant 0$ copies of $G$. For $t \geqslant 1$, \emph{$tK_1$} is called the \emph{empty graph} on $t$ vertices. For $s \geqslant 1$, the \emph{$s$-subdivision $G^{\langle s\rangle}$} of $G$ is the graph obtained by replacing every edge of $G$ with an $s$-path. So $G^{\langle 1\rangle} = G$. Let $e$ be an edge of a tree $T$ with end vertices $u$ and $v$. Let \emph{$T_e^u$} be the component of $T - e$ containing $u$. Let \emph{$T_u^e := T_e^v \cup \{e\}$}. A \emph{unit} is a vertex or an edge. The subgraph of $G$ induced by $U \subseteq V(G)$ is the maximal subgraph of $G$ with vertex set $U$. For $\emptyset \neq F \subseteq E(G)$, the subgraph of $G$ induced by $F \cup U$ is the minimal subgraph of $G$ with edge set $F$, and vertex set including $U$.

An \emph{$\ell$-arc} (or \emph{$*$-arc} if we ignore the length) is a sequence $\vec{L} := (v_0, e_1, \ldots, e_\ell, v_\ell)$, where $e_i$ is an edge with end vertices $v_{i - 1}$ and $v_i$ such that $e_{j} \neq e_{j + 1}$ for $i \in [\ell] := \{1, 2, \ldots, \ell\}$ and $j \in [\ell - 1]$. Note that $\vec{L}$ is different from $-\vec{L} := (v_{\ell}, e_\ell, \ldots, e_1, v_0)$ unless $\ell = 0$. For each $i \in [\ell]$, $\vec{e}_i := (v_{i - 1}, e_i, v_i)$ is called an \emph{arc} for short. $v_0$, $v_\ell$, $\vec{e}_1$ and $\vec{e}_\ell$ are the \emph{tail vertex}, \emph{head vertex}, \emph{tail arc} and \emph{head arc} of $\vec{L}$ respectively. The $\ell$-link (or $*$-link if we ignore the length) $L := [v_0, e_1, \ldots, e_\ell, v_\ell] = [v_{\ell}, e_\ell, \ldots, e_1, v_0]$ is obtained by taking $\vec{L}$ and $-\vec{L}$ as a single object; that is, $L := \{\vec{L}, -\vec{L}\}$. For example, a $0$-link is a vertex, and a $1$-link can be identified with an edge. For $0 \leqslant i \leqslant j \leqslant \ell$, $\vec{R} := \vec{L}(i, j) := (v_i, e_{i + 1}, \ldots, e_j, v_j)$ is called a \emph{$(j - i)$-arc} (or a \emph{subsequence} for short) of $\vec{L}$, and $\vec{L}[i, j] := R$ is a \emph{$(j - i)$-link} (or a \emph{subsequence} for short) of $L$. For $\ell \geqslant 1$, we say $L$ is \emph{formed} by the $(\ell - 1)$-links $\vec{L}[0, \ell - 1]$ and $\vec{L}[1, \ell]$. An \emph{$\ell$-dipath} is an $\ell$-arc without repeated vertices. We say $\vec{L}$ is an \emph{$\ell$-dicycle} if $v_0 = v_\ell$ and $\vec{L}(0, \ell - 1)$ is a dipath. And in this case, $L$ is an \emph{$\ell$-cycle}. For $\ell \geqslant 2$, we usually use \emph{$\vec{C}_\ell$} to denote an $\ell$-dicycle, and use \emph{$C_\ell$} to denote an $\ell$-cycle.
An \emph{$\ell$-path} is an $\ell$-link without repeated vertices. We use \emph{$\vec{\mathscr{L}}_\ell(G)$}, \emph{$\mathscr{L}_\ell(G)$}, and \emph{$\mathscr{P}_\ell(G)$} to denote the sets of $\ell$-arcs, $\ell$-links, and $\ell$-paths of $G$ respectively.

Let $\vec{L} := (v_0, e_1, \ldots, e_\ell, v_\ell) \in \vec{\mathscr{L}}_\ell(G)$, and $\vec{R} := (u_0, f_1, \ldots, f_s, u_s) \in \vec{\mathscr{L}}_s(G)$ such that $v_\ell = u_0$ and $e_\ell \neq f_1$. The \emph{conjunctions} of $\vec{L}$ and $\vec{R}$ are $\vec{Q} := (\vec{L}.\vec{R}) := (v_0, e_1, \ldots, e_\ell, v_\ell = u_0, f_1, \ldots, f_s, u_s) \in \vec{\mathscr{L}}_{\ell + s}(G)$ and $[\vec{L}.\vec{R}] := Q \in \mathscr{L}_{\ell + s}(G)$. For $\vec{Q} \in \vec{\mathscr{L}}_{\ell + s}(G)$, let $\vec{L}_i := \vec{Q}(i, \ell + i)$, and $\vec{Q}_j := \vec{Q}(j - 1, \ell + j)$ for $i \in [0, s]$ and $j \in [s]$. By definition $Q_j \in \mathscr{L}_{\ell + 1}(G)$ yields an edge $Q_j^{[\ell]} := [L_{j - 1}, Q_j, L_j]$ of $\mathbb{L}_\ell(G)$. So $Q^{[\ell]} := [L_0, Q_1, L_1, \ldots, L_{s - 1}, Q_s, L_s]$ can be seen as an $s$-link, while $\vec{Q}^{[\ell]} := (L_0, Q_1, L_1, \ldots, L_{s - 1}, Q_s, L_s)$ is an $s$-arc of $\mathbb{L}_\ell(G)$. We say that $L_0$ can be shunted to $L_s$ through $\vec{Q}$. $Q^{\{\ell\}} := \{L_0, L_1, \ldots, L_s\}$ and $\vec{Q}^{\{\ell\}} := \{\vec{L}_0, \vec{L}_1, \ldots, \vec{L}_s\}$ are the sets of \emph{images} of $L_0$ and $\vec{L}_0$ respectively during this shunting. More generally, for $R, R' \in \mathscr{L}_\ell(G)$, we say $R$ can be \emph{shunted} to $R'$ if there are $\ell$-links $R = R_0, R_1, \ldots, R_s = R'$, and $*$-arcs $\vec{P}_1, \ldots, \vec{P}_s$ of $G$ such that $R_{i - 1}$ can be shunted to $R_i$ through $\vec{P}_i$ for $i \in [s]$.

\section{Examples and basis}\label{sec_example}
We begin with some examples and basic analysis that help to build some general impressions on $\ell$-roots, and explain some of our motivations.

%For a given graph $H$, in connection with deciding if $H$ is an $\ell$-link graph, we should read, as much as possible, from $H$ the structure information of its possible minimal $\ell$-roots. In this section, we bound the number of cyclic and acyclic components, and the maximum degree of $G \in \mathbb{R}_\ell(H)$. The results will be used to bound the order of $G$ in the next section.

%\subsection{Components}
%We exemplify that, for $G \in \mathbb{R}_\ell(H)$, $c(G)$ may be different from $c(H)$.
First of all, we characterise the minimal $\ell$-roots of $2K_1$.
\begin{lem}\label{lem_barK2}
Let $P := [v_0, \ldots, v_\ell]$ be an $\ell$-path, and $T_i$ be obtained from $P$ by pasting $v_i$ at an end vertex of another $i$-path, where $i \in [0, \ell]$. Then $\mathbb{R}_\ell(2K_1) = \{2P, T_i | 1 \leqslant i \leqslant \lfloor\frac{\ell - 1}{2}\rfloor\}$. Further, $|\mathbb{R}_\ell(2K_1)|$ is $1$ if $\ell = 0$, and is $\lfloor\frac{\ell + 1}{2}\rfloor$ if $\ell \geqslant 1$.
\end{lem}
\pf Clearly, for $1 \leqslant i \leqslant \lfloor\frac{\ell - 1}{2}\rfloor$, $\mathbb{L}_\ell(2P) \cong \mathbb{L}_\ell(T_i) \cong 2K_1$. If $G \in \mathbb{R}_\ell(2K_1)$ contains a cycle $C$, then $\mathbb{L}_\ell(G)$ contains a cycle $\mathbb{L}_\ell(C)$, which is impossible. Thus $G$ is a forest containing exactly one $\ell$-path $Q$ other than $P$. If $P$ and $Q$ are vertex disjoint, then $G = P \cup Q \cong 2P$ because of the minimality. Otherwise, assign directions such that $\vec{P} = (\vec{P}_1.\vec{R}.\vec{P}_2)$ and $\vec{Q} = (\vec{Q}_1.\vec{R}.\vec{Q}_2)$, where $R$ is the maximal common path of $P$ and $Q$, $P_i \in \mathscr{P}_{s_i}(G)$ and $Q_i \in \mathscr{P}_{t_i}(G)$ for $i \in \{1, 2\}$. Since $P \neq Q$, without loss of generality, $\vec{P}_1 \neq \vec{Q}_1$ and $s_1 \geqslant t_1$. Then $s_2 \leqslant t_2$, and $\vec{L} := (\vec{P}_1.\vec{R}.\vec{Q}_2) \in \vec{\mathscr{L}}_{\ell + t_2 - s_2}(G)\setminus\{\vec{Q}\}$. Since $2K_1$ contains no edge, $\mathscr{L}_{\ell + 1}(G) = \emptyset$ and so $s_2 = t_2$. Thus $s_1 = t_1 \geqslant 1$, and $L \in \mathscr{P}_\ell(G)$. So $\vec{L} = \vec{P}$ since otherwise, $G$ contains three pairwise different $\ell$-paths $L, P$ and $Q$. Note that $[\vec{P}_1.-\vec{Q}_1] \in \mathscr{P}_{2s_1}(G)\setminus \{P, Q\}$. So $2s_1 < \ell$ and the lemma follows.
\qed

%\pf Let $t := \lfloor\frac{\ell - 1}{2}\rfloor$. Clearly, for $i \in [t]$, $2P, T_i \in \mathbb{R}_\ell(2K_1)$. Now we consider the reverse.
%If $G \in \mathbb{R}_\ell(2K_1)$ contains a cycle $C$, then $\mathbb{L}_\ell(G)$ contains a cycle $\mathbb{L}_\ell(C)$, contradicting $\mathbb{L}_\ell(G) = 2K_1$. Thus $G$ is a forest containing exactly one $\ell$-path $Q := [u_0, \ldots, u_\ell]$ other than $P$. Since $G$ is acyclic, at least one end, say $v_0$, of $P$ does not belong to $Q$. Similarly, $u_0 \notin V(P)$. If $P$ and $Q$ are disjoint, then $G = P \cup Q \cong 2P$ since $G$ is minimal. Otherwise, there are $1 \leqslant j \leqslant i \leqslant \ell$ such that $v_i = u_j$. If $v_\ell \neq u_\ell$, then $1 \leqslant j \leqslant i \leqslant \ell - 1$, and
%
%there is another path of length $\geqslant \ell$ from $v_0$ to $u_\ell$, a contradiction. So there exists $i \in [\ell - 1]$ such that $v_j \neq u_j$ iff $j < i$. Note there is a $2i$-path $\notin \{P, Q\}$ from $v_0$ to $u_0$. Thus $2i < \ell$ and the lemma follows.
%\qed

%By Whitney \cite{Whitney1932}, $\mathbb{R}_1(K_3) = \{K_3, K_{1, 3}\}$. As a generalisation, Broersma and Hoede \cite{BH1989} pointed out that a $6$-cycle $C$, $C, K_{1, 3}^{\langle 2\rangle} \in \mathbb{R}_2(C)\}$. We continue this research in the area of $\ell$-link graphs and characterise the $\ell$-roots of all cycles. Note for a given $\ell \geqslant 0$, every cycle has a unique cyclic $\ell$-root; that is, the cycle itself. So we only need to consider $\ell$-roots that are trees.

\begin{figure}[!h]
\centering
\includegraphics{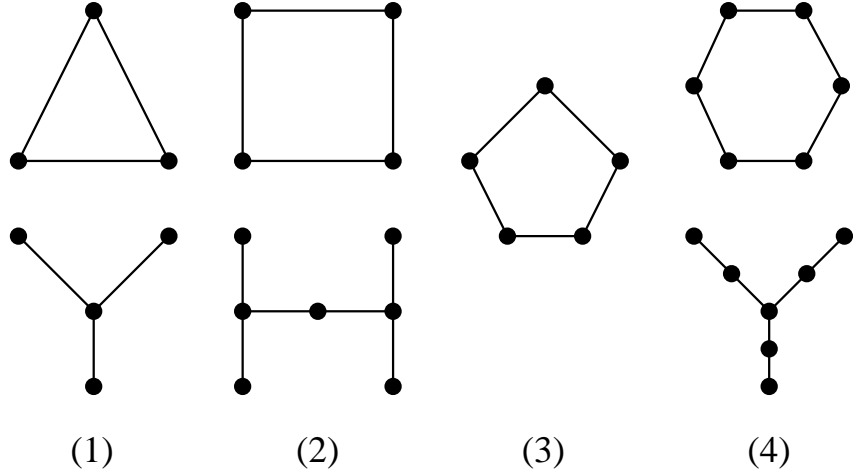}
\caption[Minimal roots of cycles]{(1) $\mathbb{R}_1(C_3)$ \quad (2) $\mathbb{R}_3(C_4)$ \quad (3) $\mathbb{R}_\ell(C_5)$ \quad (4) $\mathbb{R}_2(C_6)$}
\label{F:rootsO}
\end{figure}

Whitney \cite{Whitney1932} proved that $\mathbb{R}_1(K_3) = \{K_3, K_{1, 3}\}$ (Figure \ref{F:rootsO}(1)). As a generalisation, Broersma and Hoede \cite{BH1989} pointed out that a $6$-cycle is the $2$-path (and hence $2$-link) graph of $K_{1, 3}^{\langle 2\rangle}$ and itself (Figure \ref{F:rootsO}(4)). Figure \ref{F:rootsO}(2) are the minimal $3$-roots of $C_4$. Figure \ref{F:rootsO}(3) is the minimal $\ell$-root of $C_5$, which is $C_5$ itself. More generally, we now characterise the minimal $\ell$-roots of all cycles. Clearly, for a given $\ell \geqslant 0$, every cycle has a unique cyclic minimal $\ell$-root which is isomorphic to itself. So we only need to consider acyclic minimal $\ell$-roots.

%(as what will be explained later, a cycle may have two or more cyclic $\ell$-path roots)

%\begin{lem}
%Let $\ell \geqslant 1$, $C_\ell$ be an $\ell$-cycle, $T \in \mathbb{R}_\ell(C_t)$ be a tree, and $X := K_{3, 3} \setminus E(C_4)$. Then either $t = 3\ell$ and $T \cong K_{1, 3}^{\langle\ell\rangle}$, or $\ell$ is even, $t = 2\ell$, and $T \cong X^{\langle\ell/2\rangle}$.
%\end{lem}

%from $P := [v_0, \ldots, v_{\ell + s}]$ by pasting an end vertex of an $s$-path to $v_s$, and that of another $s$-path to $v_\ell$.

\begin{lem}\label{lem_cyclerootT}
Let $T$ be a minimal acyclic $\ell$-root of a $t$-cycle. Then $\ell \geqslant 1$, and either $t = 3\ell$ and $T \cong K_{1, 3}^{\langle\ell\rangle}$, or there is $s \geqslant 1$ such that $t = 4s$, $\ell \geqslant 2s + 1$, and $T$ is obtained by joining the middle vertices of two $2s$-paths by an $(\ell - s)$-path.
\end{lem}

Lemma \ref{lem_cyclerootT} is proved in Section \ref{sec_boundSize}. Together with Lemma \ref{lem_conClass}, it gives all $\ell$-roots $G$ of a $t$-cycle as follows: If $G$ is cyclic, it is the disjoint union of a $t$-cycle and zero or more trees of diameter at most $\ell - 1$. Otherwise, $G$ is a forest and $\ell \geqslant 1$. In this situation, if $t = 3\ell$, since $\diam(K_{1, 3}^{\langle\ell\rangle}) = 2\ell$, then $G$ is the disjoint union of $K_{1, 3}^{\langle\ell\rangle}$ and trees of diameter at most $\ell - 1$. In the final case of Lemma \ref{lem_cyclerootT}, $\diam(T) = \ell + s \leqslant 2\ell - 2$. Let $[v_0, \ldots, v_{\ell - s}]$ be the path of $T$ between the middle vertices of the two $2s$-paths. Then $\ecc_{T}(v_i) = \max\{i, \ell - i\}$ for $i \in [\ell - s]$. So $G$ is obtained from $T$ by first pasting to each $v_i$, where $i \in \{2, 3, \ldots, \ell - s - 2\}$, the root of a rooted tree of height less than $\min\{i, \ell - s - i\}$, and then adding acyclic components of diameter less than $\ell$.

In comparison with Lemma \ref{lem_cyclerootT}, a $4s$-cycle has at least three minimal $(2s + 1)$-path roots, two of which are cyclic.
Let $s \geqslant 1$ and $\ell \geqslant s + 1$ be integers. Let $G$ be the graph formed by connecting two $(s + 1)$-cycles with an $(\ell - s)$-path. One can easily check that $G$ is a minimal $\ell$-path root of a $4s$-cycle.

Broersma and Hoede \cite{BH1989} asked that, for $\ell = 2$, whether there exist three pairwise non-isomorphic simple connected graphs whose $\ell$-path graphs are isomorphic to the same connected nonnull graph. A negative answer was given by Li and Liu \cite{LiLiu2008}. We now give examples of graphs $H$ for which $|\mathbb{R}_3(H)| \geqslant 3$, and for each $\ell \geqslant 4$, there are graphs $H$ such that $|\mathbb{R}_3(H)|$ is unbounded. The following construction will be useful.

\begin{figure}[!h]
\centering
\includegraphics{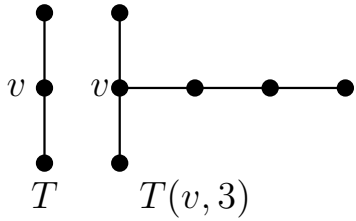}
\caption[A tree is an $\ell$-link graph for large enough $\ell$]{The $3$-link graph of $T(v, 3)$ is $T$}
\label{F:rootTreeShunting}
\end{figure}

Let $T$ be a finite tree, and $v$ be a vertex of degree $d$ in $T$. Assign to $v$ an integer $t_v$ as follows: If $d \geqslant 2$, then $t_v := \diam(T)$. If $d \leqslant 1$ and $T$ is a path, then $t_v := -1$. Otherwise, $d = 1$, and there exists a path $[v, \ldots, e, u]$ of minimum length such that $\deg_T(u) \geqslant 3$. In this case, let $t_v := \diam(T_e^u)$. Denote by $T(v, \ell)$ the tree obtained by pasting an end vertex of an extra $\ell$-path to $v$. For example, in Figure \ref{F:rootTreeShunting}, $T$ is a $2$-path with middle vertex $v$, and $T(v, 3)$ is the tree obtained from $T$ by pasting an end vertex of a $3$-path at $v$. It is not difficult to see that the $3$-link graph $T(v, 3)$ is isomorphic to $T$.

Figure \ref{F:rootTreeShunting} is just a special case of the following lemma which says that every tree $T$ is an $\ell$-link graph for each $\ell > \diam(T)$.

\begin{lem}\label{lem_copylemma}
Let $T$ be a tree, and $\ell \geqslant t_v + 1$. Then $T \cong \mathbb{L}_\ell(T(v, \ell)) = \mathbb{P}_\ell(T(v, \ell))$.
\end{lem}
\pf Let $\vec{L}$ be the $\ell$-arc of head vertex $v$ such that $L$ is the extra path. Consider the shunting of $L$ in $G := T(v, \ell)$. One can check that the mapping $\vec{L}'[\ell, \ell] \mapsto L'$, for every image $\vec{L}'$ of $\vec{L}$, is an isomorphism from $T$ to $\mathbb{L}_\ell(G)$.
\qed

An \emph{orbit} of a graph $G$ is a maximal subset $U \subseteq V(G)$ such that for every pair of vertices in $U$, one can be mapped to the other by an automorphism of $G$ (see \cite{Biggs1993} for more about algebraic graph theory). For each $\ell \geqslant 1$, the number of non-isomorphic trees $T(v, \ell)$, over all $v \in V(T)$, equals the number of orbits of $T$. Two vertices in the same orbit have the same eccentricity. So the number of orbits of $T$ is at least the radius plus one, which is $\lfloor \diam(T)/2\rfloor + 1$, with equality holds if and only if the set of leaves is an orbit of $T$.

As explained below, Lemma \ref{lem_copylemma} implies that there are infinitely many trees $T$ of diameter $3$ such that $\mathbb{Q}_\ell(T) \cap \mathbb{R}_\ell(T)$ contains at least four trees, where $\ell \geqslant 3$.

Let $T$ be obtained by adding an edge between the centers of $K_{1, p}$ and $K_{1, q}$, where $p > q \geqslant 1$. Then $\diam(T) = 3$, and $T$ has four orbits. So there are exactly four non-isomorphic $T(v, \ell)$'s for each $\ell \geqslant 3$.

Let $\ell \geqslant s \geqslant 4$ and $k \geqslant 1$ be given integers. Lemma \ref{lem_copylemma} also implies that there exists a tree $T$ of diameter $s$, such that $\mathbb{Q}_\ell(T) \cap \mathbb{R}_\ell(T)$ contains at least $k$ trees: Let $T(k)$ be obtained by pasting a leaf of each star $K_{1, i + 1}$, where $i \in [k]$, at the same end vertex of an $(s - 2)$-path. Then $\diam(T(s)) = s$, and the number of orbits of $T(k)$ is $\lfloor\frac{s}{2}\rfloor + 1$ if $k = 1$, and $s + 2k - 1$ if $k \geqslant 2$.

\section{Constructing $\ell$-equivalence classes}\label{sec_consEquiClass}
In this section, we explain the process of constructing $\ell$-roots from minimal $\ell$-roots, which allows us to concentrate on the latter in our future study.

\subsection{Incidence units}
Two $*$-links of a graph $G$ are \emph{incident} if one is a subsequence of the other. A $*$-link is said to be \emph{$\ell$-incident} if it is incident to an $\ell$-link. It follows from the definitions immediately that every $\ell$-link is $\ell$-incident, and every $\ell$-incident $*$-link is $s$-incident, for $s \leqslant \ell$. Conversely, a $t$-link is not $\ell$-incident if and only if it is not $s$-incident for each $s \geqslant \ell$. And if this is the case, then $\ell \geqslant t + 1$. Let $T$ be a tree obtained by joining the centres of two stars such that each star contains at least two edges. Then all units of $T$ are $3$-incident. However, every $2$-path of a star is not $3$-incident in $T$. The fact below allows us to focus on incidence units of trees of finite diameter.

%\begin{exa}\label{exa_vertexroots}
%Let $Q := [u_0, \ldots, u_4]$ be a path, and $T_0$ be obtained by adding an extra neighbour $v_0$ to $u_1$. For $i \geqslant 1$, let $T_i$ be obtained from $T_0$ by pasting the center of $K_{1, i}$ with $u_2$. Then for each $i \geqslant 0$, $T_i$ is connected, and $\mathbb{L}_4(T_i) \cong 2K_1$. Further $Q$ and all units of the path $P := [u_0, u_1, v_0]$ are $4$-incident in $T_i$. However, neither $P$ or any leaf or edge of $K_{1, i}$ is $4$-incident in $T_i$.
%\end{exa}

\begin{lem}\label{lem_cycleIncident}
Let $G$ be a connected nonnull graph. Then $G$ contains a cycle or a ray if and only if for every $\ell \geqslant 0$, all units of $G$ are $\ell$-incident.
\end{lem}
\pf
$(\Leftarrow)$ Suppose not. Then $G$ is a tree of finite diameter $s$. Then no unit of $G$ is $(s + 1)$-incident. $(\Rightarrow)$ Let $X$ be a cycle or a ray in $G$. Clearly, every unit of $X$ is $\ell$-incident. So we only need to show that every $e \in E(G)\setminus E(X)$ is $\ell$-incident. Since $G$ is connected, there exists a dipath $\vec{P}$ of minimum length with tail edge $e$ and head vertex $x \in V(X)$. Clearly, $X$ contains an $\ell$-arc $\vec{R}$ starting from $x$. Then $L := (\vec{P}.\vec{R})[0, \ell]$ is an $\ell$-link of $G$ incident to $e$.
\qed

%%$G$ contains a cycle or a ray iff for every $\ell \geqslant 0$, each unit of $G$ is $\ell$-incident.
%\pf
%$(\Leftarrow)$ Suppose not. Then $G$ is a tree of finite diameter. Hence as long as $\ell > \diam(G)$, no unit of $G$ is $\ell$-incident. $(\Rightarrow)$ Let $X$ be a cycle or a ray of $G$. For each vertex $v \in V(G)$, let $\vec{P} := \vec{P}_G(X, v)$ be a dipath of minimum length $s$ from $v$ to $X$. Let $\vec{L} := \vec{L}_\ell(X, v)$ be the initial subsequence of length $\ell$ of an $(\ell + s)$-arc of $X \cup P$ starting from $\vec{P}$ and end at at $X$. Then $L_\ell(X, v) := L$ is an $\ell$-link of $G$ incident to $v$. For each arc $\vec{e}$ of $G$ with head $v$, if $e$ is on $\vec{L}$ or $\ell \leqslant 1$, there is nothing to show. Otherwise, $[\vec{e}.\vec{L}(0, \ell - 1)]$ is an $\ell$-link incident to $e$.
%\qed

Let $T$ be a tree. By definitions, if $\diam(T) \leqslant \ell - 1$, then none unit of $T$ is $\ell$-incident.
The following statement says that, to study $\ell$-incidence units of $T$, we only need to consider the case of $\ell \leqslant \diam(T) \leqslant 2\ell - 4$.
\begin{lem}\label{lem_dlr4}
Let $T$ be a tree of diameter at least $\max\{\ell, 2\ell - 3\}$. Then all units of $T$ are $\ell$-incident.
\end{lem}
\pf
Let $s := \diam(T)$, and $Q$ be an $s$-path of $T$. Since $s \geqslant \ell$, every unit of $Q$ is $\ell$-incident. So we only need to show that each $e \in E(T)\setminus E(Q)$ is $\ell$-incident in $T$. Let $\vec{P}$ be a dipath of minimum length $t$ with tail edge $e$ and head vertex some $u \in V(Q)$. Then $t \geqslant 1$, and $P$ and $Q$ are edge disjoint. Clearly, $Q$ has a dipath $\vec{R}$ of length at least $\lceil s/2 \rceil \geqslant \ell - 1$ from $u$. So $e$ is incident to the path $[\vec{P}.\vec{R}]$ of length at least $t + \ell - 1 \geqslant \ell$.
\qed

%Let $G$ be a connected nonnull graph and $\ell \geqslant 0$. Then $G$ contains a cycle or $\diam(G) \geqslant \max\{\ell, 2\ell - 3\}$ if and only if all units of $G$ are $\ell$-incident.

%$s$ and radius $r$. Then $\ell \geqslant s + 1$ if and only if none unit of $T$ is $\ell$-incident. If $\ell \leqslant r$ or $\ell \leqslant \min\{3, s\}$, then all units of $T$ are $\ell$-incident.

%Let $T$ be a tree, and $[u, e, v]$ be a $1$-link of $T$. Denote by $T_u^e$ the subtree induced by $T_e^v$ and $e$.

Wu et al. \cite{WuChao2004} presented a linear time algorithm computing the eccentricity of a vertex of a finite tree. Based on this work, the following observation provides a linear time algorithm testing if a vertex is $\ell$-incident in a finite tree.
\begin{obs}\label{obs_vertexincidence}
Let $T$ be a tree and $\ell \geqslant 0$ be an integer. Then $u \in V(T)$ is $\ell$-incident in $T$ if and only if either $u$ is a leaf and $\ecc_T(u) \geqslant \ell$, or there exist different $e, f \in E(T)$ incident to $u$, such that $\ecc_{T_u^e}(u) + \ecc_{T_u^f}(u) \geqslant \ell$.
\end{obs}

Based on Observation \ref{obs_vertexincidence}, the lemma below can be formalised into a linear time algorithm for testing if an edge of a finite tree is $\ell$-incident.
\begin{lem}\label{lem_elcovered}
Let $\ell \geqslant 0$, and $P$ be a path of a tree $T$. Then all units of $P$ are $\ell$-incident in $T$ if and only if both end vertices of $P$ are $\ell$-incident in $T$.
\end{lem}
\pf
We only need to consider $(\Leftarrow)$ with the length of $P$ at least $1$. The case of $\ell \leqslant 3$ follows from Lemma \ref{lem_dlr4}. Now let $\ell \geqslant 4$. For a contradiction, let $P$ be a minimal counterexample such that its ends $u, v$ are contained in two  $\ell$-paths $Q_u$ and $Q_v$ respectively. Clearly, $Q_u$ contains a sub-path $L_u$ starting from $u$ and of length $s_u \geqslant \lceil\ell/2\rceil$. By the minimality of $P$, none inner vertex of $P$ belongs to $Q_u$ or $Q_v$. So the union of $L_u, P$ and $L_v$ forms a path of length at least $s_u + s_v + 1 > \ell$ in $T$, contradicting that $P$ is not $\ell$-incident.
\qed
%\noindent {\bf Remark.} As exemplified in Example \ref{exa_vertexroots}, $P$ may not be $\ell$-incident even when all its units are $\ell$-incident.

\subsection{Incidence subgraphs}
The \emph{$\ell$-incidence subgraph $G[\ell]$} of a graph $G$ is the graph induced by the $\ell$-incidence units of $G$. By definition, $G[\ell] = G$ if $\ell = 0$ or $G$ is null. And for each $\ell \in \{1, 2, 3\}$, $G[\ell]$ can be obtained from $G$ by deleting all acyclic components of diameter less than $\ell$. For each $G \in \mathbb{R}_\ell[K_1]$, $G[\ell]$ is an $\ell$-path. The statements below, follow from the definitions and Lemma \ref{lem_cycleIncident}, allow us to concentrate on incidence subgraphs of trees of finite diameter.
\begin{coro}\label{coro_roughGl}
Let $G$ be a graph, and $s \geqslant \ell \geqslant 0$ be integers. Then every $s$-link of $G$ belongs to $G[\ell]$. If further $G$ is nonnull and connected, then $G$ contains a cycle or a ray if and only if for every $\ell \geqslant 0$, $G = G[\ell]$.
\end{coro}

%As a consequence of Lemma \ref{lem_dlr4}, the following statements suggest us to focus on the situation that $\diam(T) \geqslant \ell \geqslant \max\{4, \radi(T)\}$.
%\begin{coro}\label{coro_dlr4}
%Let $T$ be a tree, $\ell \geqslant 0$ be an integer. Then $\ell \geqslant \diam(T) + 1$ iff $T[\ell]$ is null. Otherwise, if $\ell \leqslant \radi(T)$, or $\ell \leqslant \min\{3, \diam(T)\}$, then $T = T[\ell]$.
%\end{coro}
%
%By the analysis above, when study $\ell$-incident subgraphs, we can focus on a tree $T$ of diameter $\geqslant \ell \geqslant \max\{4, \radi(T) + 1\}$.

A rough structure of $T[\ell]$ can be derived from Lemma \ref{lem_elcovered}.
\begin{coro}\label{coro_Tltree}
Let $T$ be a tree of finite diameter, and $\ell \leqslant \diam(T)$. Then $T[\ell]$ is an induced subtree of $T$. And each leaf of $T[\ell]$ is a leaf of $T$.
\end{coro}
\pf
By Lemma \ref{lem_elcovered}, $T[\ell] \leqslant T$. Let $v$ be a leaf of $T[\ell]$. By Corollary \ref{coro_roughGl}, there is an $\ell$-path $L$ of $T[\ell]$ with an end $v$. Suppose for a contradiction that $v$ is not a leaf of $T$. Then there exists $e \in E(T) \setminus E(T[\ell])$ incident to $v$. Then the union of $e$ and $L$ forms an $(\ell + 1)$-path of $T$, contradicting that $e$ is not $\ell$-incident.
\qed

Let $u \in V(T)$ and $X$ be a subtree of $T$. Denote by $T_{X}^u$ the component of $T - E(X)$ containing $u$. An accurate structure of $T[\ell]$ is given as follows.
%Corollary \ref{coro_Tltree} indicates that $u \in V(T[\ell])$ is the only unit of $T_{T[\ell]}^u$ that is $\ell$-incident in $T$.
\begin{lem}\label{lem_TlTu}
Let $T$ be a tree of finite diameter, $\ell \leqslant \diam(T)$, and $X$ be a subtree of $T$. Then $X = T[\ell]$ if and only if $X = X[\ell]$, and for each $u \in V(X)$, either
\begin{itemize}
\item[{\bf (1)}] $\ecc_{X}(u) \geqslant \ell - 1$, and $T^u := T_{X}^u$ is a single vertex $u$. Or
\item[{\bf (2)}] $\lceil\ell/2\rceil \leqslant \ecc_{X}(u) \leqslant \ell - 2$, and $\ecc_{X}(u) + \ecc_{T^u}(u) \leqslant \ell - 1$.
\end{itemize}
\end{lem}
\pf If $\ell \leqslant 3$, then by Lemma \ref{lem_dlr4}, $T = T[\ell]$ and the lemma follows. Now let $\ell \geqslant 4$. $(\Rightarrow)$ By Corollary \ref{coro_roughGl}, we have $s := \diam(X) = \diam(T) \geqslant \ell$, and every $u \in V(X)$ is $\ell$-incident in $X$. So $X[\ell] = X$ is nonnull, and $s \geqslant \ecc_{X}(u) \geqslant \radi(X) = \lceil s/2 \rceil \geqslant \lceil\ell/2 \rceil$. By Corollary \ref{coro_Tltree}, $T^u$ is a maximal subtree of $T$, of which the only unit that is $\ell$-incident in $T$ is the vertex $u$. Let $v \in V(T^u)$ such that $t := \dist(u, v) = \ecc_{T^u}(u)$. Then $T^u$ is not a single vertex if and only if $t \geqslant 1$. If this is the case, then $\ell - 1 \geqslant \ecc_T(v) \geqslant t + \ecc_{X}(u)$, and the statement follows.

$(\Leftarrow)$ Since $X = X[\ell] \subseteq T[\ell]$, we only need to show that $T[\ell] \subseteq X$. Suppose not. Then there exists some $\vec{P} := (v_0, e_1, \ldots, e_\ell, v_\ell) \in \vec{\mathscr{L}}_\ell(T)$, and a maximum $s \in [\ell]$ such that $\vec{P}[0, s]$ belongs to $T^{u}$, where $u := v_s$.

According to {\bf (2)}, $\radi(T^{u}) \leqslant \ecc_{T^{u}}({u}) \leqslant \ell - 1 - \ecc_{X}(u) \leqslant \ell - 1 - \lceil\ell/2\rceil = \lfloor\ell/2\rfloor - 1$. So $\diam(T^{u}) \leqslant 2\radi(T^{u}) \leqslant \ell - 2$. Hence there is a maximum $t \geqslant s + 1$ such that $t \leqslant \ell$ and $\vec{R} := \vec{P}[s, t]$ belongs to $X$. Since $\ecc_{X}(u) + \ecc_{T^{u}}(u) \leqslant \ell - 1$, we have $t \leqslant \ell - 1$. Since $X = X[\ell]$ is nonnull, there exists $\vec{L} := (\vec{L}_1.\vec{P}(s_1, t_1).\vec{L}_2) \in \vec{\mathscr{L}}_\ell(X)$, where $s \leqslant s_1 < t_1 \leqslant t$, and $L_1$ and $L_2$ are edge disjoint with $P$.

Since $\ecc_{X}(u) + \ecc_{T^{u}}(u) < \ell$, $(\vec{P}(0, t_1).\vec{L}_2)$ is a dipath of length less than $\ell$. So $\vec{L}_2$ is of length less than $\ell - t_1$. Since $\ecc_{X}(v_t) + \ecc_{T^{v_t}}(v_t) < \ell$, we have $(\vec{L}_1.\vec{P}(s_1, \ell))$ is a dipath of length less than $\ell$. Thus $\vec{L}_1$ is of length less than $s_1$. But then $\vec{L}$ is of length less than $s_1 + t_1 - s_1 + \ell - t_1 = \ell$, a contradiction.
\qed

\subsection{Equivalence classes}
In this subsection we build the relationships among $\ell$-minimal graphs, $\ell$-incidence subgraphs, and $\ell$-equivalence classes.
\begin{lem}\label{lem_minEq}
Each $\ell$-finite graph $G$ is $\ell$-equivalent to $G[\ell]$.
\end{lem}
\pf
By Corollary \ref{coro_roughGl}, all $\ell$-links and $(\ell + 1)$-links of $G$ belong to $G[\ell]$. So $\mathbb{L}_\ell(G) = \mathbb{L}_\ell(G[\ell])$. Note that $G[\ell] \subseteq G$. Thus the lemma follows.
\qed

The following lemma links $\ell$-incidence units with $\ell$-minimal graphs.
\begin{lem}\label{lem_minInci}
An $\ell$-finite graph $G$ is $\ell$-minimal if and only if $G = G[\ell]$.
\end{lem}
\pf
By Corollary \ref{coro_roughGl}, every unit of $G[\ell]$ is $\ell$-incident in $G[\ell]$. So deleting a unit from $G[\ell]$ erases at least one $\ell$-link from $G[\ell]$. So $G[\ell]$ is $\ell$-minimal. Conversely, if $G[\ell] \subset G$, then $G$ is not $\ell$-minimal since, by Lemma \ref{lem_minEq}, $\mathbb{L}_\ell(G[\ell]) = \mathbb{L}_\ell(G)$.
\qed

Below we connect $\ell$-equivalence relation and $\ell$-incidence graphs.
\begin{lem}\label{lem_minUniq}
Given $\ell$-finite graphs $X, Y$, $X \sim_\ell Y$ if and only if $X[\ell] \cong Y[\ell]$.
\end{lem}
\pf
$(\Leftarrow)$ Let $Z := X[\ell] \subseteq X, Y$. By Lemma \ref{lem_minEq}, $\mathbb{L}_\ell(X) \cong \mathbb{L}_\ell(Z) \cong \mathbb{L}_\ell(Y)$. So $X \sim_\ell Y$. $(\Rightarrow)$ By definitions, there exists an $\ell$-minimal graph $Z \subseteq X, Y$ such that $\mathbb{L}_\ell(X) \cong \mathbb{L}_\ell(Y) \cong \mathbb{L}_\ell(Z)$. By Lemma \ref{lem_minInci}, $Z = Z[\ell] \subseteq X[\ell]$ since $Z \subseteq X$. But by Lemma~\ref{lem_minEq}, $\mathbb{L}_\ell(Z) = \mathbb{L}_\ell(X) = \mathbb{L}_\ell(X[\ell])$. So $Z \cong X[\ell]$ since, by Lemma~\ref{lem_minInci}, $X[\ell]$ is $\ell$-minimal. Similarly, $Y[\ell] \cong Z$ and the lemma follows.
\qed

We show in the following that $\sim_\ell$ is an equivalence relation.

\medskip

\noindent {\bf Proof of Lemma \ref{lem_miniUniq}.} The reflexivity and symmetry of $\sim_\ell$ follow from the definition. To show the transitivity, let $X \sim_\ell Y$ and $Y \sim_\ell Z$. Then by Lemma~\ref{lem_minUniq}, $X[\ell] \cong Y[\ell] \cong Z[\ell]$, and hence $X \sim_\ell Z$. The uniqueness of the $\ell$-minimal graph in its class follows from Lemmas \ref{lem_minInci} and \ref{lem_minUniq}. The fact that $G[\ell]$ is an induced subgraph of $G$ follows from Corollaries \ref{coro_roughGl} and \ref{coro_Tltree}. \qed

Now we construct $\ell$-roots from minimal $\ell$-roots.

\medskip

\noindent {\bf Proof of Lemma \ref{lem_conClass}.} Let $Z \sim_\ell G$, and $Y$ be a component of $Z$. Since $Z, Y$ and $G$ are $\ell$-finite, they do not contain rays. If $Y$ contains a cycle, then $Y = Y[\ell]$ is a component of $G$ by Corollary \ref{coro_roughGl}. Now let $Y$ be a tree. If $\diam(Y) \leqslant \ell - 1$, then $Y[\ell] \cong K_0$. In this case, there can be arbitrarily many such $Y$. If $\diam(Y) \geqslant \max\{\ell, 2\ell - 3\}$, then $Y = Y[\ell]$ is a component of $G$ by Lemma \ref{lem_dlr4}.  The case of $\ell \leqslant \diam(Y) \leqslant 2\ell - 4$ follows from Lemma \ref{lem_TlTu}.\qed

\section{Partitioned $\ell$-link graphs}
A sufficient and necessary condition for an $\ell$-link graph to be connected was given by Jia and Wood \cite{JiaWood2013}. In this section, we study the cyclic components of $\ell$-roots. The investigation helps to further understand the structure of $\ell$-link graphs, and bound the parameters of minimal $\ell$-roots.

%As a , we bound the maximum degree of minimal $\ell$-roots.
%\subsection{Partitioned graphs}

\subsection{Definitions and basis} Let $H$ be a graph admitting a partition $\mathcal{V}$ of $V(G)$ and a partition $\mathcal{E}$ of $E(G)$. Then $\tilde{H} := (H, \mathcal{V}, \mathcal{E})$ is called a \emph{partitioned graph}. For each graph $G$, let \emph{$\mathcal{V}_0(G) := \{\{v\} \subseteq V(G)\}$}, and \emph{$\mathcal{E}_0(G) := \{\{e\} \subseteq E(G)\}$}. Let $\ell \geqslant 1$. For $R \in \mathscr{L}_{\ell - 1}(G)$, let \emph{$\mathscr{L}_{\ell + 1}(R)$} be the set of $\ell$-links of $G$ of middle subsequence $R$, $\mathscr{L}^{[\ell]}(R) := \{Q^{[\ell]}| Q \in \mathscr{L}_{\ell + 1}(R)\}$, and \emph{$\mathcal{E}_\ell(G) := \{\mathscr{L}^{[\ell]}(R) \neq \emptyset | R \in \mathscr{L}_{\ell - 1}(G)\}$}. Let $E_G(u, v)$ be the set of edges of $G$ between $u, v \in V(G)$, and \emph{$\mathcal{V}_1(G) := \{E_G(u, v) \neq \emptyset | u, v \in V(G)\}$}. For $\ell \geqslant 2$, let \emph{$\mathcal{V}_\ell(G) := \{\mathscr{L}_\ell(R) \neq \emptyset| R \in \mathscr{L}_{\ell - 2}(G)\}$}. By \cite[Lemma 4.1]{JiaWood2013}, for $\ell \neq 1$, $\mathcal{V}_\ell(G)$ consists of independent sets of $\mathbb{L}_\ell(G)$. For $\ell \geqslant 0$, $\widetilde{\mathbb{L}}_\ell(G) := (\mathbb{L}_\ell(G), \mathcal{V}_\ell(G), \mathcal{E}_\ell(G))$ is a partitioned graph, and called a \emph{partitioned $\ell$-link graph} of $G$. $(\mathcal{V}_\ell(G), \mathcal{E}_\ell(G))$ is called an \emph{$\ell$-link partition} of $H$. $G$ is an \emph{$\ell$-root} of $\tilde{H}$ if $\widetilde{\mathbb{L}}_\ell(G) \cong \tilde{H}$. Denote by \emph{$\mathbb{R}_\ell[\tilde{H}]$} and \emph{$\mathbb{R}_\ell(\tilde{H})$} respectively the sets of all  $\ell$-roots and minimal $\ell$-roots of $\tilde{H}$. The statement below follows from definitions.

%Two partitioned graphs $\tilde{H}$ and $\tilde{X}$ are said to be \emph{isomorphic} if there is a graph isomorphism $\varphi: H \mapsto X$, such that for every pair of units $u, v$ of $H$, $u, v$ are in the same part of $\tilde{H}$ iff $\varphi(u)$ and $\varphi(v)$ are in the same part of $\tilde{X}$. If there exists some graph $G$ such that $\tilde{H} \cong \widetilde{\mathbb{L}}_\ell(G)$, then $\tilde{H}$ is said to be a partitioned $\ell$-link graph (of $G$). And in this case, $(\mathcal{V}, \mathcal{E})$ is called an \emph{$\ell$-link partition} of $H$, while $G$ is called an \emph{$\ell$-root} of $\tilde{H}$.

%; that is, it can be organised into a partitioned $\ell$-link graph
%It is meaningful to know that the number \emph{$p_\ell(H)$} of nonisomorphic $\ell$-link partitions of $H$ is finite.

\begin{propo}
A graph is an $\ell$-link graph if and only if it admits an $\ell$-link partition. Moreover, $\mathbb{R}_\ell(H)$ (respectively, $\mathbb{R}_\ell[H]$) is the union of $\mathbb{R}_\ell(\tilde{H})$ (respectively, $\mathbb{R}_\ell[\tilde{H}]$) over all partitioned graphs $\tilde{H}$ of $H$.
\end{propo}

An \emph{$\ell$-link} (respectively, \emph{$\ell$-arc}) of $\tilde{H}$ is an $\ell$-link (respectively, $\ell$-arc) of $H$ whose consecutive edges are in different edge parts of $\tilde{H}$. The lemma below indicates that every $s$-link of $\widetilde{\mathbb{L}}_\ell(G)$ arises from an $(\ell + s)$-link of $G$.

\begin{lem}\label{lem_sHslG}
Let $G$ be a graph, and $L$ be an $s$-link of $\mathbb{L}_\ell(G)$. Then $L$ is an $s$-link of $\widetilde{\mathbb{L}}_\ell(G)$ if and only if there exists an $(\ell + s)$-link $R$ of $G$ such that $L = R^{[\ell]}$.
\end{lem}
\pf
It is trivial for $\ell = 0$ or $s \leqslant 1$. Now let $\ell \geqslant 1$, $s \geqslant 2$, and $L := [L_0, Q_1, \ldots, Q_s, L_s]$. $(\Leftarrow)$ Let $L = R^{[\ell]}$, and $P_i$ be the middle $(\ell - 1)$-link of $Q_i$, where $i \in [s]$. By definition $Q_i^{[\ell]} = [L_{i - 1}, Q_i, L_i] \in \mathscr{L}^{[\ell]}(P_i) \in \mathcal {E}_\ell(G)$. By \cite[Observation 3.3]{JiaWood2013}, $P_i \neq P_{i + 1}$ for $i \in [s - 1]$. So $Q_{i}^{[\ell]}$ and $Q_{i + 1}^{[\ell]}$ are in different edge parts of $\mathcal {E}_\ell(G)$. Thus $L$ is an $s$-link of $\widetilde{\mathbb{L}}_\ell(G)$.

$(\Rightarrow)$ Let $L$ be an $s$-link of $\widetilde{\mathbb{L}}_\ell(G)$, where $L_i := [v_i, e_{i + 1}, \ldots, v_{i + \ell}] \in \mathscr{L}_\ell(G)$ for $i \in \{0, 1\}$, and $Q_1 := [v_0, e_1, \ldots, e_{\ell + 1}, v_{\ell + 1}] \in \mathscr{L}_{\ell + 1}(G)$. Suppose $Q_2$ has the form $[u_0, f_1, v_1, \ldots, e_{\ell + 1}, v_{\ell + 1}]$. Then $Q_2$ and $Q_1$ have the same middle $(\ell - 1)$-link $P := [v_1, e_2, \ldots, e_\ell, v_\ell]$, and hence are in the same part $\mathscr{L}^{[\ell]}(P) \in \mathcal {E}_\ell(G)$, contradicting that $L$ is an $s$-link of $\widetilde{\mathbb{L}}_\ell(G)$. So $Q_2$ has the form $[v_1, e_2, \ldots, e_{\ell + 2}, v_{\ell + 2}]$, and $L_2 = [v_2, e_3, \ldots, e_{\ell + 2}, v_{\ell + 2}]$. Continue this analysis, we have that $L_i = [v_i, e_{i + 1}, \ldots, v_{i + \ell}] \in \mathscr{L}_\ell(G)$, and $Q_i = [v_{i - 1}, e_i, \ldots, v_{i + \ell}] \in \mathscr{L}_{\ell + 1}(G)$ for $i \in [s]$. Then $R := [v_0, e_1, \ldots, v_{s + \ell}] \in \mathscr{L}_{\ell + s}(G)$ and  $L = R^{[\ell]}$.
\qed

As shown below, each closed $s$-link of $\widetilde{\mathbb{L}}_\ell(G)$ stems from a closed $s$-link of $G$.

%\begin{lem}\label{lem_cycleGH}
%Let $\ell \geqslant 0$, $s \geqslant 2$, $G$ be a graph, and $[L_0, Q_1, \ldots, Q_s, L_s = L_0]$ be a closed $s$-link of $\tilde{H} := \widetilde{\mathbb{L}}_\ell(G)$. Then there exists an $(\ell + s)$-link $[v_0, e_1, \ldots, e_{\ell + s}, v_{\ell + s}]$ such that $v_i = v_{s + i}$ for $i \in \{0, \ldots, \ell\}$, $e_i = e_{s + i}$ for $i \in [\ell]$, $L_i = [v_i, e_{i + 1}, \ldots, e_{i + \ell}, v_{i + \ell}]$ for $i \in \{0, \ldots, s\}$, and $Q_i = [v_{i - 1}, e_i, \ldots, e_{i + \ell}, v_{i + \ell}]$ for $i \in [s]$.
%\end{lem}

%Let $L$ and $\vec{R}$ be described in Lemma \ref{lem_sHslG}, and $s \geqslant 2$. Then $L$ is closed iff $\vec{R}(0, \ell) = \vec{R}(s, \ell + s)$.

\begin{lem}\label{lem_cycleGH}
Let $\ell \geqslant 0$ and $s \geqslant 2$ be integers. Let $G$ be a graph, and $\vec{R}$ be an $(\ell + s)$-arc of $G$. Then $R^{[\ell]}$ is a closed $s$-link of $\widetilde{\mathbb{L}}_\ell(G)$ if and only if $\vec{R}(0, \ell) = \vec{R}(s, \ell + s)$.
\end{lem}
\pf By Lemma \ref{lem_sHslG}, $[L_0, Q_1, \ldots, Q_s, L_s] := R^{[\ell]}$ is an $s$-link of $\widetilde{\mathbb{L}}_\ell(G)$. Clearly, $R^{[\ell]}$ is closed if and only if $\vec{R}[0, \ell] = \vec{R}[s, \ell + s]$. Let $\vec{R} := (v_0, e_1, \ldots, v_{\ell + s})$. Suppose $\vec{R}(0, \ell) = \vec{R}(\ell + s, s)$; that is, $(v_0, e_1, \ldots, e_\ell, v_\ell) = (v_{\ell + s}, e_{\ell + s}, \ldots, e_{s + 1}, v_s)$. Then $Q_1 = [v_0, e_1, \ldots, e_{\ell + 1}, v_{\ell + 1}]$ and $Q_s = [v_{s - 1}, e_{s}, \ldots, e_{\ell + s}, v_{\ell + s}]$ have the same middle $(\ell - 1)$-link $P := [v_1, e_{2}, \ldots, e_\ell, v_\ell] = [v_{\ell + s - 1}, e_{\ell + s - 1}$, $\ldots, e_{s + 1}, v_s]$. So $Q_1, Q_s$ correspond to edges of $\widetilde{\mathbb{L}}_\ell(G)$ from the same edge part $\mathscr{L}^{[\ell]}(P)$, contradicting that $R^{[\ell]}$ is an $s$-link of $\widetilde{\mathbb{L}}_\ell(G)$. Thus $\vec{R}(0, \ell) = \vec{R}(s, \ell + s)$.
\qed

%A \emph{cycle} of $\tilde{H}$ is a cycle of $H$ of which the consecutive edges are in different parts of $\mathcal {E}$. As a common sense, $H$ contains a closed $*$-link if and only if it contains a cycle. But this may not be true for $\tilde{H}$.
%\begin{exa}
%For $i \in [2]$, let $\vec{C}_i := (v_0^i, e_1^i, \ldots, e_5^i, v_5^i = v_0^i)$ be a cycle. Let $H$ be the \emph{bicycle graph} (see, for example, \cite{Zaslavsky1998}) obtained by joining $C_1$ and $C_2$ with a common vertex $v := v_0^1 = v_0^2$. Let $V := \{v\}$, and $E_3^i := \{e_3^i\}$. For $j \in [2]$, let $V_j^i := \{v_j^i, v_{5 - j}^i\}$, and $E_j^i := \{e_j^i, e_{6 - j}^i\}$. Let $\mathcal {V} := \{V, V_j^i | i, j \in [2]\}$, and $\mathcal {E} := \{E_j^i | i \in [2], j \in [3]\}$. Then $[\vec{C}_1.\vec{C}_2]$ is a closed $10$-link of $\tilde{H} := (H, \mathcal {V}, \mathcal {E})$. However, $\tilde{H}$ contains no cycle.
%\end{exa}

\subsection{Cycles in partitioned graphs}
A \emph{cycle} of $\tilde{H}$ is a cycle of $H$ whose the consecutive edges are in different edge parts. $\tilde{H}$ and its partition are \emph{cyclic} if $\tilde{H}$ contains a cycle, and \emph{acyclic} otherwise. For example, for each $t$-cycle $C$, $\widetilde{\mathbb{L}}_\ell(C)$ is cyclic. When $t \geqslant 3$ is divisible by $3$ or $4$, by Lemma \ref{lem_cyclerootT}, there exists a tree $T$ and an integer $\ell$ such that $C \cong \mathbb{L}_\ell(T)$ can be organised into an acyclic partitioned graph $\widetilde{\mathbb{L}}_\ell(T)$. Each component $X$ of $H$ corresponds to a partitioned subgraph $\tilde{X}$ of $\tilde{H}$. $\tilde{X}$ is called a \emph{component} of $\tilde{H}$. Let \emph{$o(\tilde{H})$} and \emph{$a(\tilde{H})$} be the cardinalities of the sets of cyclic and acyclic components of $\tilde{H}$ respectively. The following lemma says that the number of cyclic components is invariant (see, for example, \cite{Olver1999}) under the partitioned $\ell$-link graph construction.

\begin{lem}\label{lem_oHeqoG}
Let $G$ be a graph, and $\ell \geqslant 0$. Then $o(G) = o(\widetilde{\mathbb{L}}_\ell(G))$.
\end{lem}
\pf
Let $\tilde{H} := \widetilde{\mathbb{L}}_\ell(G)$, $X$ be a component of $G$ containing a cycle $C$, $\tilde{X}_\ell$ be the component of $\tilde{H}$ containing $\widetilde{\mathbb{L}}_\ell(C)$. We only need to show that $\varphi: X \mapsto \tilde{X}_\ell$ is a bijection from the cyclic components of $G$ to that of $\tilde{H}$. For $i \in \{1, 2\}$, let $C_i$ be a closed $s_i$-link of $\tilde{H}$. By Lemma \ref{lem_cycleGH}, $G$ contains an $(\ell + s_i)$-arc $\vec{R}_i$ such that $C_i = R_i^{[\ell]}$, and that $\vec{R}_i(0, s_i)$ is a closed $s_i$-arc $\vec{O}_i$ of $G$.

First we show that $\varphi$ is a well defined surjection. Assume that $R_1$ and $R_2$ are $*$-links of $X$. Note that $\vec{R}_i[0, \ell]$ can be shunted into $O_i$ of $X$. Thus by \cite[Lemma 3.7]{JiaWood2013}, $\vec{R}_1[0, \ell]$ can be shunted to $\vec{R}_2[0, \ell]$ in $X$, and the images of the shunting form a $*$-link from $C_1$ to $C_2$ in $\tilde{X}_\ell$.

We still need to show that $\varphi$ is injective. Assume $C_1$ and $C_2$ are joined by a $*$-link $Q$ of $X_\ell$ between, say, $\vec{R}_1[0, \ell] \in V(C_1)$ and $\vec{R}_2[0, \ell] \in V(C_2)$. Then $\vec{R}_i[0, \ell] \in \mathscr{L}_\ell(G)$, where $i \in \{1, 2\}$, can be shunted to each other in $G$, with images corresponding to the vertices of $Q$ in $X_\ell$. \qed

Every graph is a disjoint union of its connected components. But the relationship between a partitioned graph $\tilde{H}$ and its components is more complicated. For different components $\tilde{X}$ and $\tilde{Y}$ of $\tilde{H}$, it is possible that a vertex part $U$ of $\tilde{X}$ and a vertex part $V$ of $\tilde{Y}$ are two disjoint subsets of a vertex part $W$ of $\tilde{H}$. Lemma \ref{lem_oHeqoG} leads to a rough process of building $\widetilde{\mathbb{L}}_\ell(G)$ from its components. We can first fix all cyclic components such that no two units from different components are in the same part of $\tilde{H}$. And then for each acyclic component, we either set it as an independent fixed component, or merge some of its vertex parts with that of a unique fixed graph to get a larger fixed graph.

%Certainly, there are many different choices in the construction process, which lead to different partitioned graphs. $\tilde{H}$ is just one of them. And some resulted partitioned graphs may not be partitioned $\ell$-link graphs.

%A very useful corollary can be derived from Lemma \ref{lem_oHeqoG} immediately:
%\begin{coro}\label{coro_olGeq0Gtree}
%An $\ell$-finite graph $G$ is acyclic if and only if $o(\tilde{\mathbb{L}}_\ell(G)) = 0$.
%\end{coro}

By definition, $o(\tilde{H}) \leqslant o(H)$, $a(H) \leqslant a(\tilde{H})$, and $a(H) + o(H) = c(H) = c(\tilde{H}) = a(\tilde{H}) + o(\tilde{H})$. As a consequence of Lemma \ref{lem_oHeqoG}, we have:
\begin{coro}\label{coro_aHeqcH} For each integer $\ell \geqslant 0$, a graph $G$ is acyclic if and only if
$a(\widetilde{\mathbb{L}}_\ell(G)) = c(\mathbb{L}_\ell(G))$.
\end{coro}

\subsection{Computing cyclic components}
We explain how to decide whether a component of $\tilde{H} := (H, \mathcal{V}, \mathcal{E})$ is cyclic, and compute $o(\tilde{H})$ and $a(\tilde{H})$ in quadratic time.
Let \emph{$\mathcal{E}(u)$} be the edge parts in $\mathcal{E}$ incident to $u \in V(H)$, and \emph{$r({\mathcal{E}}) := \max\{|\mathcal{E}(u)|\;| u \in V(H)\}$}.

%Clearly, $|\mathcal{E}(u)| \leqslant \deg_H(u)$.

%By \cite[Lemma 4.1]{JiaWood2013}, $|\mathcal{E}(u)| \leqslant 2$ if $\tilde{H}$ is a partitioned $\ell$-link graph with $\ell $.

\begin{defi}\label{defi_link2arcgraph}
Let $\tilde{H} := (H, \mathcal{V}, \mathcal{E})$ be a partitioned graph, and $\vec{H}$ be the digraph with vertices $(u, E)$, where $E \in \mathcal{E}(u)$, such that there is an arc from $(u, E)$ to $(v, F)$ if $E \neq F$, and there is $e \in E$ between $u$ and $v$.
\end{defi}

% A graph homomorphism $\varphi: \vec{H} \mapsto H$ is $\varphi(u, E) = u$, and $\varphi((u, E), (v, F)) = e$

%\begin{defi}\label{defi_link2arcgraph}
%Let $\tilde{H} := (H, \mathcal{V}, \mathcal{E})$ be a partitioned graph. For each $u \in V(H)$, if $k := |\mathcal{E}(u)| \leqslant 1$, then add $2 - k$ pendent edges $e$ that are incident to $u$, and put each $\{e\}$ into $\mathcal{E}$. Note $\mathcal{E}$ is changed, but $H$ is not. Now for each $u \in V(H)$, $|\mathcal{E}(u)| = 2$. Let $\vec{H}$ be the digraph with vertices $(u, E)$, where $u \in V(H)$ and $E \in \mathcal{E}(u)$, such that there is an arc from $(u, E)$ to $(v, F)$ if $E \neq F$, and there is an edge of $E$ between $u$ and $v$. A graph homomorphism $\varphi: \vec{H} \mapsto H$ is $\varphi(u, E) = u$, and $\varphi((u, E), (v, F)) = e$, where $u, v \in V(H)$, $\mathcal{E}(v) = \{E, F\}$, and $e \in E$ between $u$ and $v$.
%\end{defi}

%\begin{lem}\label{lem_p2acyclic}
%In Definition \ref{defi_link2arcgraph}, let $\tilde{X}$ be a component of $\tilde{H}$. Then $\vec{X} := \varphi^{-1}(\tilde{X})$ is the disjoint union of a set of components of $\vec{H}$. Moreover, $\tilde{X}$ is cyclic if and only if $\vec{X}$ contains a directed cycle.
%\end{lem}

In the following, we transfer the problem of computing $o(\tilde{H})$ to that of detecting if a component of $\vec{H}$ contains a dicycle.

\begin{lem}\label{lem_p2acyclic}
In Definition \ref{defi_link2arcgraph}, for each component $\tilde{X}$ of $\tilde{H}$, we have that $\vec{X}$ is the disjoint union of some components of $\vec{H}$. And $\vec{H}$ is the disjoint union of $\vec{X}$ over all components $\tilde{X}$ of $\tilde{H}$.  Moreover, $\tilde{X}$ is cyclic if and only if $\vec{X}$ contains a dicycle.
\end{lem}
\pf The first two statements follow from definitions.
We now prove the last statement.
Let $s \geqslant 2$, and $C := [v_0, e_1, \ldots, e_s, v_s = v_0]$ be a cycle of $X$. For $i \in [s]$, let $e_i \in E_i \in \mathcal{E}$. Let $e_{s + 1} := e_1$ and $E_{s + 1} := E_1$. By definition, $C$ is a cycle of $\tilde{X}$ if and only if $E_i \neq E_{i +1}$ for $i \in [s]$, which is equivalent to saying that $u_{i - 1} := (v_{i - 1}, E_i)$ is a vertex of $\vec{H}$, while $\vec{f}_i := (u_{i - 1}, u_i)$ is an arc of $\vec{H}$ for $i \in [s + 1]$; that is, $(u_0, \vec{f}_1, \ldots, \vec{f}_s, u_s = u_0)$ is a dicycle of $\vec{X}$.
\qed
\noindent {\bf Remark.} Let $n := n(H)$, $m := m(H)$ and $r := r(\mathcal{E})$. Then we have $n(\vec{H}) = \sum_{u \in V(H)}|\mathcal{E}(u)| \leqslant \sum_{u \in V(H)}\deg_H(u) = 2m$. Every arc $(u, e, v)$ of $H$ with $e \in E \in \mathcal{E}$ corresponds to $|\mathcal{E}(v)| - 1$ arcs of $\vec{H}$; that is, $((u, E), (v, F))$ for $F \in \mathcal{E}(v)\setminus \{E\}$. So $m(\vec{H}) \leqslant \sum_{(u, e, v)\in \vec{\mathscr{L}}_1(H)}(|\mathcal{E}(v)| - 1) \leqslant \sum_{v \in V(H)}\deg_H(v)(|\mathcal{E}(v)| - 1) \leqslant 2m(r - 1)$. An $O(n + m)$-time algorithm for dividing $H$ into connected components was given by Hopcroft and Tarjan~\cite{HopcroftTarjan1973}. For each component $\tilde{X}$ of $\tilde{H}$, Tarjan's algorithm~\cite{Tarjan1972}, with time complexity $O(n(\vec{X}) + m(\vec{X})) = O(rm(X))$, can be used to detect the existence of dicycles in $\vec{X}$, and hence that of cycles in $\tilde{X}$ by Lemma \ref{lem_p2acyclic}. So the time complexity for computing $o(\tilde{H})$ is $O(m^2 + n)$ in general situations, and is $O(m + n)$ if $r$ is bounded. By \cite[Lemma 4.1]{JiaWood2013}, $r(\mathcal{E}_\ell(G)) \leqslant 2$ for $\ell \geqslant 1$. So it requires $O(m + n)$-time to compute the number of cyclic components of a finite partitioned $\ell$-link graph for each $\ell \geqslant 0$.

\section{Bounding the number of minimal roots}\label{sec_boundSize}
We bound in this section the order, size, maximum degree and total number of minimal $\ell$-roots of a finite graph. This lays a basis for solving the recognition and determination problems for $\ell$-link graphs in our future work.

%The results are the basis for solving the recognition and determination problems for $\ell$-link graphs \cite{JiaWood2013rec}.

%As a byproduct, we show that trees (finite or infinite) of bounded diameters are better-quasi-ordered under induced subgraph relation.

\subsection{Incidence pairs}
For $s, \ell \geqslant 0$, and $L \in \mathscr{L}_\ell(G)$, let \emph{$\mathcal {I}_G(L, s)$} be the set of pairs $(L, R)$ such that $R \in \mathscr{L}_s(G)$ is incident to $L$. Let \emph{$i_G(L, s) := |\mathcal {I}_G(L, s)|$}. Define \emph{$\girth(L)$} to be $+\infty$ if $L$ is a path, and the minimum length of a sub cycle of $L$ otherwise. To dodge confusions, denote by \emph{$\hat{L}$} the graph induced by the units of $L$. Then $\girth(\hat{L}) \leqslant \girth(L)$ and the inequality may hold. For example, let $L := [v_0, e_1, v_1, e_2, v_2, e_3, v_0, e_4, v_1]$. Then $\girth(\hat{L}) = 2$ and $\girth(L) = 3$.
\begin{lem}\label{lem_inciCounting}
Let $\ell \geqslant s \geqslant 0$, $L \in \mathscr{L}_\ell(G)$, and $g := \girth(L^{[s]})$. Then $\min\{g, \ell - s + 1\} \leqslant i_G(L, s) \leqslant \ell - s + 1$. Further, $i_G(L, s) = \ell - s + 1$ if and only if $g = +\infty$ if and only if $L^{[s]}$ is an $(\ell - s)$-path if and only if $g \geqslant \ell - s + 1$. Otherwise, if $g \leqslant \ell - s$, then $i_G(L, s) = g$ if and only if $\hat{L^{[s]}}$ is a $g$-cycle.
\end{lem}
\pf
$L^{[s]}$ is an $(\ell - s)$-link of $H := \mathbb{L}_s(G)$. So $t := i_G(L, s) =  |L^{\{s\}}| \leqslant \ell - s + 1$, with the last equality holds if and only if $L^{[s]}$ is a path. If $g \leqslant \ell - s$, then $L^{[s]}$ contains a sub $g$-cycle on which every pair of different vertices of $H$ corresponding to a pair of different sub $s$-links of $L$. So $t \geqslant g$, with equality holds if and only if all units of $H$ on $L^{[s]}$ belong to the $g$-cycle.
\qed

% In particular, $i_G(L, \ell - 1) = 2$.

Let \emph{$\mathcal {I}_G(\ell, s) := \bigcup_{L \in \mathscr{L}_\ell(G)}\mathcal {I}_G(L, s)$}, and \emph{$i_G(\ell, s):= |\mathcal {I}_G(\ell, s)|$}. Then $i_G(\ell, s) = \sum_{L \in \mathscr {L}_\ell(G)}$ $i_G(L, s)$ can be bounded as follows:

%In particular, $i_G(\ell, \ell - 1) = 2|\mathscr{L}_\ell(G)|$.

\begin{coro}\label{coro_inciCounting}
Let $\ell \geqslant s \geqslant 0$ be integers, and $G$ be an $\ell$-finite graph of girth $g$. Then
\begin{eqnarray*}\label{eq_gllink}
\min\{g, \ell - s + 1\}|\mathscr{L}_\ell(G)| &\leqslant& i_G(\ell, s) \leqslant (\ell - s + 1)|\mathscr{L}_\ell(G)|.
\end{eqnarray*}
Further, $i_G(\ell, s) = (\ell - s + 1)|\mathscr{L}_\ell(G)|$ if and only if $g \geqslant \ell - s + 1$. If $g \leqslant \ell - s$, then $i_G(\ell, s) = g|\mathscr{L}_\ell(G)|$ if and only if $G$ is a disjoint union of $g$-cycles.
\end{coro}
\pf
Let $t: = \min\{\girth(L^{[s]})| L\in \mathscr{L}_\ell(G)\}$. Note that $g \geqslant \ell - s + 1$ if and only if $t = +\infty$. And $g \leqslant \ell - s$ if and only if $t = g$. So the statements can be verified by summing the results in Lemma \ref{lem_inciCounting} over $L \in \mathscr{L}_\ell(G)$.
\qed

The order and size of minimal $\ell$-roots are bounded as follows.
\begin{lem}\label{lem_HboundG}
Let $\ell \geqslant 1$ and $G$ be an $\ell$-minimal graph. Then both $G$ and $\tilde{H} := \widetilde{\mathbb{L}}_\ell(G)$ are finite. Further, $m(G) \leqslant \ell n(H)$, and $n(G) \leqslant \ell n(H) + a(\tilde{H})$.
\end{lem}
\pf
Since $G$ is $\ell$-finite, $\tilde{H}$ is finite. By Corollary \ref{coro_inciCounting}, $i_G(\ell, 1) \leqslant \ell n(H)$ is finite. By Lemma \ref{lem_minInci}, for each $e \in E(G)$, $i_G(e, \ell) \geqslant 1$. Summing this inequation over $e \in E(G)$, we have $m(G) \leqslant i_G(\ell, 1)$ is finite. By Lemma \ref{lem_oHeqoG}, $o(G) = o(\tilde{H})$, and hence $a(G) = c(G) - o(G) \leqslant c(\tilde{H}) - o(\tilde{H}) = a(\tilde{H})$.
So $n(G) \leqslant m(G) + a(G) \leqslant \ell n(H) + a(\tilde{H})$ is finite.
\qed

\subsection{The number of minimal roots} We have known that $|\mathbb{R}_0(G)| = 1$, and $|\mathbb{R}_\ell(K_s)| = 1$ for $s \in \{0, 1, 2\}$. By Lemma \ref{lem_barK2}, $|\mathbb{R}_\ell(2K_1)| = \lfloor\frac{\ell + 1}{2}\rfloor$ for $\ell \geqslant 1$. Let $\tilde{a}_\ell(H)$ be the maximum $a(\tilde{H})$ over all partitioned $\ell$-link graphs $\tilde{H} := (H, \mathcal{V}, \mathcal{E})$. Clearly, $\tilde{a}_\ell(H) \leqslant c(H)$. Denote by \emph{$\psi(p, q)$} the number of nonisomorphic graphs with $p$ vertices and $q$ edges. Then $\psi(p, q)$ is $0$ if $p \leqslant 1$ and $q \geqslant 1$; is $1$ if $p = 2$ or $q = 0$; is $1$ if $p \geqslant 2$ and $q \leqslant 1$; is less than $\binom{p}{2}^q$ if $p \geqslant 3$ and $q \geqslant 1$; and by Hardy~\cite{Hardy1920}, is the nearest integer to $\frac{(q + 3)^2}{12}$ if $p = 3$.
Further, when $p \geqslant 3$ and $q \geqslant 1$, we have that $\binom{p}{2} \geqslant 3$, and $\binom{p - 1}{2}/\binom{p}{2} = (1 - \frac{2}{p})^{\frac{p}{2}\cdot\frac{2}{p}} < \textbf{e}^{\frac{-2}{p}}$. Thus $\sum_{i = 0}^p\sum_{j = 0}^q \psi(i, j) \leqslant p + q + \sum_{i = 3}^p\sum_{j = 1}^q\binom{i}{2}^j < p + q + \frac{3}{2}\sum_{i = 3}^p\binom{i}{2}^q < p + q + \frac{3}{2}\binom{p}{2}^q[1 - \textbf{e}^{\frac{-2q}{p}}]^{-1}$.

\begin{lem}\label{lem_lminiFinite}
Let $\ell \geqslant 1$ be an integer, and $H$ be a finite graph. Let $n := n(H) \geqslant 2$, and $a := \tilde{a}_\ell(H)$. Then $H$ has at most $(\ell n + a)^{2\ell n}$ minimal $\ell$-roots, of which at most $\frac{3}{2}(\ell n + 1)^{\ell n - 1}$ are trees, and at most $\frac{3}{4}a(a + 1)(\ell n + a)^{\ell n - 1}$ are forests.
\end{lem}
\pf Let $G \in \mathbb{R}_\ell(H)$. By Lemma \ref{lem_HboundG}, $m(G) \leqslant \ell n$, and $n(G) \leqslant \ell n + a \leqslant 2\ell n$. By the analysis above, $|\mathbb{R}_\ell(H)| < 2\ell n + a + \frac{3}{2}\binom{\ell n + a}{2}^{\ell n}(1 - \textbf{e}^{\frac{-2\ell n}{\ell n + a}})^{-1} \leqslant 2\ell n + a + \frac{3\textbf{e}}{2(\textbf{e} - 1)}\binom{\ell n + a}{2}^{\ell n} < (\ell n + a)^{2\ell n}$.

%Let $t := {n(\delta - 1)^{1 - \ell}}$. By Lemma \ref{lem_boundingsizeroots}(3), the number of graphs in $\mathbb{R}_\ell(H)$ of minimum degree at least $\delta \geqslant 2$ is at most $$\sum_{n = 1}^t\sum_{m = n}^t \binom{n}{2}^{m} \leqslant t - 1 + 2\sum_{n = 3}^t\binom{n}{2}^t \leqslant 2t\binom{t}{2}^t \leqslant 2\sqrt{2}(\frac{t}{\sqrt{2}})^{2t + 1}.$$

Cayley's formula \cite{Borchardt1860,Cayley2009} states that there are $p^{p - 2}$ unequal trees on vertex set~$[p]$. So the number of trees in $\mathbb{R}_\ell(H)$ is at most $\sum_{p = 1}^{\ell n + 1}p^{p - 2} < \frac{3}{2}(\ell n + 1)^{\ell n - 1}$.

By Aigner and Ziegler \cite{AignerZiegler2010}, the number of unequal forests on vertex set $[p]$ of $k$ components is $kp^{p - k - 1}$. So the number of forests in $\mathbb{R}_\ell(H)$ is at most $\sum_{k = 1}^a\sum_{p = a}^{\ell n + k}kp^{p - k - 1} < \frac{3}{2}\sum_{k = 1}^{a}k(\ell n + k)^{\ell n - 1} < \frac{3}{4}a(a + 1)(\ell n + a)^{\ell n - 1}$.
\qed

\subsection{The maximum degree of minimal roots} We have proved that minimal $\ell$-roots of a finite graph are finite. So in this subsection, we only deal with finite graphs. Let $\mathcal{E}$ be a partition of $E(H)$. We may identify $E \in \mathcal{E}$ with the subgraph of $H$ induced by $E$. Let \emph{$D(\mathcal {E})$} be the set of $\deg_E(v)$ over $E \in \mathcal{E}$ incident to $v \in V(H)$. Let \emph{$D(G) := \{\deg_G(v) - 1 \geqslant 1 | v \in V(G)\}$}. Clearly, $D(G) = \emptyset$ if and only if $\Delta(G) \leqslant 1$ if and only if for all $\ell \geqslant 1$, $D(\mathcal{E}_\ell(G)) = \emptyset$.

\begin{lem}\label{lem_Dsdecrease}
Let $\ell \geqslant 1$, and $G$ be a finite connected graph of $\Delta(G) \geqslant 2$. Let $D := D(G)$, and $D_\ell := D(\mathcal {E}_\ell)$. If $G$ is cyclic, then $D = D_\ell$. Otherwise, $G$ is a tree of diameter $s \geqslant 2$, and $D = D_1 \supseteq \ldots \supseteq D_{s - 1} \supset \emptyset = D_{s} = D_{s + 1} = \ldots$.
\end{lem}
\pf Let $\ell, d \geqslant 1$. By definition, $d \in D_\ell$ if and only if there is an $\ell$-arc $\vec{L}$ starting from some $v \in V(G)$ with $\deg_G(v) = d + 1$. So $D = D_1$. When $\ell \geqslant 2$, $\vec{L}(0, \ell - 1)$ is an $(\ell - 1)$-arc starting from $v$. So $d \in D_{\ell - 1}$. In another word, $D_{\ell - 1} \supseteq D_\ell$ for all $\ell \geqslant 2$. 

On the one hand, let $C$ be a cycle of $G$. For each $v \in V(G)$, since $G$ is connected, there is a dipath $\vec{P}$ of minimum length from $v$ to some $u \in V(C)$. Clearly, there is an $\ell$-arc $\vec{R}$ of $C$ starting from $u$. Then $\vec{L} := (\vec{P}.\vec{R})(0, \ell)$ is an $\ell$-arc of $G$ starting from $v$. Thus $D_\ell = D$ by the analysis above.

On the other hand, let $G$ be a tree of diameter $s \geqslant 2$. By the definitions, $\ell \geqslant s$ if and only if $E(\mathbb{L}_\ell(G)) = \emptyset$ if and only if $D_\ell = \emptyset$.
\qed

Define \emph{$\Delta(\mathcal {E})$} to be $0$ if $\mathcal {E} = \emptyset$, and the maximum of $D(\mathcal {E})$ otherwise. A lower bound of $\Delta(G)$ for $G \in \mathbb{R}_\ell[\tilde{H}]$ follows from Lemma \ref{lem_Dsdecrease}:

%\begin{coro}\label{coro_DsdecreaseCyc}
%Let $\ell \geqslant 1$, $\Delta := \Delta(G)$, and $\Delta_\ell := \Delta(\mathcal{E}_\ell(G))$. Then $\Delta = \Delta_\ell = 0$ if and only if $E(G) = \emptyset$. Conversely, $E(G) \neq \emptyset$ if and only if $\Delta \geqslant \Delta_\ell + 1$. Further assume that $G$ is connected. Then $G$ is cyclic if and only if for all $\ell \geqslant 1$,
%$\Delta = \Delta_\ell + 1 \leqslant \Delta(H)$.
%\end{coro}

\begin{coro}\label{coro_DsdecreaseCyc}
Let $\ell \geqslant 1$ be an integer, and $G$ be a finite graph. Let $\mathcal{E} := \mathcal{E}_\ell(G)$. Then $\Delta(G) = \Delta(\mathcal{E}) = 0$ if and only if $E(G) = \emptyset$. And $E(G) \neq \emptyset$ if and only if $\Delta(G) \geqslant \Delta(\mathcal{E}) + 1$. Further assume that $G$ is connected. Then $G$ contains a cycle if and only if for all $\ell \geqslant 1$, we have
$\Delta(G) = \Delta(\mathcal{E}) + 1 \leqslant \Delta(\mathbb{L}_\ell(G))$.
\end{coro}

Below we display a connection between the degrees of an $\ell$-minimal tree and the number of components of the $\ell$-link graph of the tree.

\begin{lem}\label{lem_Tdc}
Let $\ell \geqslant 1$ be an integer, $T$ be a finite $\ell$-minimal tree, and $v$ be a vertex of eccentricity less than $\ell$ in $T$. Then $\deg_T(v) \leqslant c(\mathbb{L}_\ell(T)) + 1$.
\end{lem}
\pf Since $T$ is $\ell$-minimal, so $\diam(T) \geqslant \ell \geqslant 1$. Thus $s := \ecc_T(v) \geqslant 1$ and $\ell \geqslant s + 1 \geqslant 2$. Let $d := \deg_T(v)$. If $d \leqslant 1$, there is nothing to show. 

Now let $d \geqslant 2$. For $i \in [d]$, let $\vec{e}_i := (v, e_i, u_i)$ be the arcs of $T$ starting from $v$. Then there exists $\vec{R} \in \vec{\mathscr{L}}_s(T)$ starting from, say, $\vec{e}_d$. Since $T$ is $\ell$-minimal, $e_i$ is $\ell$-incident in $T$. For $i \in [d - 1]$, by Lemma \ref{lem_elcovered}, $t_i := \ecc_{T_{e_i}^{u_i}}(u_i) \geqslant \ell - s - 1 \geqslant 0$. So there is a $t_i$-arc $\vec{Q}_i$ from $u_i$ to some $v_i$ in $T_{e_i}^{u_i}$. Obviously, $L_i := [\vec{Q}_i(\ell - s - 1, 0). \vec{e}_i. \vec{R}]$, for $i \in [d - 1]$, are $d - 1$ different $\ell$-paths containing $v$ in $T$. 

Suppose for a contradiction that $L_i$ can be shunted to $L_j$ for some $1 \leqslant i < j \leqslant d - 1$. Since $e_i$ separates $v$ from $v_i$, so $v$ is an image of $v_i$ during the shunting. But then $\ecc_T(v) \geqslant \ell$, a contradiction. So $L_i$ and $L_j$ correspond to vertices in different components of $H := \mathbb{L}_\ell(T)$. Hence $d - 1 \leqslant c(H)$.
\qed

As an application of Lemma \ref{lem_Tdc}, we characterise the minimal roots of a cycle.

\medskip

\noindent{\bf Proof of Lemma \ref{lem_cyclerootT}} Since $C := \mathbb{L}_\ell(T)$ is connected and $T$ is minimal, so $T$ is connected and hence is a tree. For $u, v \in V(T)$, define $D(u, v) := \deg_T(u) + \deg_T(v)$. Since $C$ is $2$-regular, $D(u, v) = 4$ if $\dist_T(u, v) = \ell$.

We claim that, if $\dist_T(u, v) \geqslant \ell$, then $D(u, v) \leqslant 4$. Suppose for a contradiction that $D(u, v) \geqslant 5$. Without loss of generality, assume $\deg_T(u) \geqslant 3$. Since $\dist_T(u, v) \geqslant \ell$, there exists some vertex $w$ on the path of $T$ from $u$ to $v$ such that $\dist_T(u, w) = \ell$. Then $\deg_T(w) \geqslant 2$. Thus $D(u, w) = 5 > 4$, a contradiction.

By the minimality of $T$, for each leaf $w$ of $T$, there exists some $v \in V(T)$ such that $\dist_T(v, w) = \ell$ and $\deg_T(v) = 3$.
By Lemma \ref{lem_Tdc}, for each $v \in V(T)$ with $\ecc_T(v) < \ell$, we have $\deg_T(v) \leqslant c(C) + 1 = 2$. So $\deg_T(v) \in \{1, 2, 3\}$ for $v \in V(T)$. Let $k$ be the number of degree-$3$ vertices in $T$. Then $k \geqslant 1$ since $T$ contains degree-$1$ vertices.

If $k = 1$, then $T$ contains exactly three leaves (\cite[Page 67]{ChartrandOellermann1993}). Since $T$ is minimal, each leaf $u$ of $T$ is the end vertex of some $\ell$-path $P$ of $T$. Let $v$ be the other end of $P$. Then $\deg(v) = 3$ since $\dist_T(u, v) = \ell$. So each leaf of $T$ is
 at distance $\ell$ from the unique degree-$3$ vertex $v$ of $T$. As a consequence, $T \cong K_{1, 3}^{\langle\ell\rangle}$ and $t = 3\ell$.

Now assume that $k \geqslant 2$. Let $\vec{Q} := (v_0, \ldots, v_q)$ be a dipath in $T$ such that $v_0$ and $v_q$ are degree-$3$ vertices at maximum distance in $T$. Since $D(v_0, v_q) = 6$, we have $q = \dist_T(v_0, v_q) < \ell$.

If $k = 2$, then $T$ contains exactly four leaves \cite[Page 67]{ChartrandOellermann1993}. Consequently, $T$ is the union of two paths $[\vec{L}_i.\vec{Q}.\vec{R}_i]$, where $i \in \{1, 2\}$, and $L_i, R_i$ are four internally disjoint paths of length $\ell_i, s_i \geqslant 1$ respectively. By the analysis above, there is an $\ell$-path between each leaf and one of $v_0$ and $v_q$. Let $(w_0, w_1, \ldots, w_{\ell_1} = v_0) := \vec{L}_1$. If $\dist_T(w_0, v_0) = \ell$, then $\dist_T(w_1, v_q) = \ell - 1 + q \geqslant \ell$, and $D(w_1, v_q) > 4$, a contradiction. So $\dist_T(w_0, v_q) = \ell$, and $\ell_1 = \ell - q$. Let $s := \ell - q$. By symmetry, $\ell_1 = \ell_2 = s_1 = s_2 = s$. Now let $\vec{P}_{ij} := (\vec{L}_i.\vec{Q}.\vec{R}_j)$ for $i, j \in \{1, 2\}$. Then the images during the shunting of $L$ through $\vec{P}_{11}, -\vec{P}_{21}, \vec{P}_{22}, -\vec{P}_{12}$ form a $4s$-cycle $C'$ of $\mathbb{L}_\ell(T)$. Hence $C' = C$, and all $\ell$-links of $T$ are aforementioned images of $L$. Thus $[\vec{L}_1.\vec{L}_2] \in \mathscr{P}_{2s}(T)$ contains no $\ell$-link; that is, $\ell \geqslant 2s + 1$.

Assume for a contradiction that $k \geqslant 3$. Since $q$ is maximum, there exists some $p \in [q - 1]$ such that $\deg_T(v_p) = 3$. By Lemma \ref{lem_Tdc}, $\ecc_T(v_p) \geqslant \ell$. So there exists an $\ell$-dipath $\vec{L} := (u_0, \ldots, u_\ell = v_p)$ of $T$ such that $v_0, v_q$ are separated from $u_0$ by $v_p$. Then $\dist_T(v_0, u_1) = \ell - 1 + p \geqslant \ell$, and $D(v_0, u_1) > 4$, which is a contradiction.
\qed

We now bound the maximum degree of a finite tree in terms of $\widetilde{\mathbb{L}}_\ell(T)$.

\begin{coro}\label{coro_DeltaT}
Let $\ell \geqslant 1$, $T$ be a finite tree, $\tilde{H} := (H, \mathcal {V}, \mathcal {E}) := \widetilde{\mathbb{L}}_\ell(T)$ and $s := \max\{\ecc_T(v) |$ $\deg_T(v) = \Delta(T)\}$. Then
\begin{itemize}
\item[{\bf (1)}] If $s \geqslant \ell + 1$, then $\Delta(T) = \Delta(\mathcal {E}) + 1 \leqslant \Delta(H)$.
\item[{\bf (2)}] If $\ell = s$, then either $s  = \Delta(T) = 1$, and $\Delta(H) = 0$; or $s \geqslant 2$ and $\Delta(T) = \Delta(\mathcal {E}) + 1 \leqslant \Delta(H) + 1$.
\item[{\bf (3)}] If $\ell \geqslant s + 1$ and $T$ is $\ell$-minimal, then $\Delta(T) \leqslant a(\tilde{H}) + 1 \leqslant c(H) + 1$.
\end{itemize}
\end{coro}
\pf
{\bf (1)} and {\bf (2)} are implied by the proof of Lemma \ref{lem_Dsdecrease}. {\bf (3)} follows from Lemma \ref{lem_Tdc} and Corollary \ref{coro_aHeqcH}.
\qed

%Consider {\bf (3)}. Clearly, $s \geqslant 1$ and $\ell \geqslant 2$. Since $T$ is $\ell$-minimal, $d := \Delta(T) \geqslant 2$. For $i \in [d]$, let $\vec{e}_i := (v, e_i, u_i)$ be the arcs of $T$ starting from $v$. Since $s = \ecc_T(v)$, there exists $\vec{R} \in \vec{\mathscr{L}}_s(T)$ starting from, say, $\vec{e}_d$. Since $G$ is $\ell$-minimal, $e_i$ is $\ell$-incident in $G$. So for $i \in [d - 1]$, by Lemma \ref{lem_elcovered}, $t_i := \ecc_{T_{e_i}^{u_i}}(u_i) \geqslant \ell - s - 1 \geqslant 0$. Thus there exists a $t_i$-arc $\vec{Q}_i$ of $T_{e_i}^{u_i}$ from $u_i$ to some $v_i \in V(T)$. Then $L_i := [\vec{Q}_i(\ell - s - 1, 0). \vec{e}_i. \vec{R}]$, for $i \in [d - 1]$, are $d - 1$ different $\ell$-paths containing $v$ in $G$. Suppose for a contradiction that $L_i$ can be shunted to $L_j$ for some $1 \leqslant i < j \leqslant d - 1$. Since $e_i$ separates $v$ from $v_i$, $v$ is an image of $v_i$ during the shunting. But then $\ecc_T(v) \geqslant \ell$, contradicting $\ell > s$. So $\mathbb{L}_\ell(L_i)$ and $\mathbb{L}_\ell(L_j)$ are in different components of $H$. Thus $d - 1 \leqslant c(H) = a(\tilde{H})$ by Corollary \ref{coro_aHeqcH}.

For $\tilde{H} := (H, \mathcal{V}, \mathcal{E})$, let $b(\tilde{H}) := \max\{a(\tilde{H}), \Delta(\mathcal{E})\}$. By Corollaries \ref{coro_DsdecreaseCyc} and \ref{coro_DeltaT}, we have that $b(\tilde{H}) \leqslant \max\{c(H), \Delta(H)\}$. Below we bound the maximum degree of minimal $\ell$-roots $G$ of $\tilde{H}$ or $H$. The results in turn help to bound the parameters of $s$-link graphs of $G$. Clearly, for $s \geqslant 1$, a graph of $q \geqslant 2$ edges contains at most $q(q - 1)^{s - 1}$ $s$-links, with equality holds if and only if all these edges are between the same pair of vertices.
\begin{lem}
Let $\ell$ and $s$ be positive integers, and $H$ be a finite graph. Let $n := n(H) \geqslant 2$, and $b := b(\tilde{H})$. Let $G \in \mathbb{R}_\ell(\tilde{H})$, and $(X, \mathscr{V}, \mathscr{E}) := \widetilde{\mathbb{L}}_s(G)$. Then $\Delta(G) \leqslant b + 1$, $\Delta(\mathscr{E}) \leqslant b$, $\Delta(X) \leqslant 2b$, $\max\{|V|, |E|\, |\, V \in \mathscr{V}, E \in \mathscr{E}\} \leqslant b^2$, $m(X) \leqslant \ell n(\ell n - 1)^s$, and $n(X) \leqslant \ell n(\ell n - 1)^{s - 1}$.
\end{lem}
\pf
Let $Y, Z$ be the subgraphs induced by the cyclic and acyclic components of $G$ respectively. Let $(H, \mathcal{V}, \mathcal{E}) := \tilde{H}$. By Corollary \ref{coro_DsdecreaseCyc}, we have $\Delta(Y) \leqslant \Delta(\mathcal{E}) + 1$. By Corollary \ref{coro_DeltaT}, we have $\Delta(Z) \leqslant a(\tilde{H}) + 1$. So $\Delta(G) = \max\{\Delta(Y), \Delta(Z)\} \leqslant b + 1$. The rest of the lemma follows from the analysis above.
\qed

\section{Path graphs}\label{sec_pathGraphs}
Some ideas and techniques used in the investigation of $\ell$-link graphs can be applied to the study of $\ell$-path graphs. We end this paper by bounding the parameters of minimal $\ell$-path roots.
\subsection{Quantitative analysis}
A graph is \emph{$\ell$-path finite} if its $\ell$-path graph is finite. Note that an $\ell$-finite graph is $\ell$-path finite, but not vise versa. For example, the disjoint union of infinitely many $\ell$-cycles is $\ell$-path finite but not $\ell$-finite.
Two $\ell$-path finite graphs $X$ and $Y$ are \emph{$\ell$-path equivalent}, write $X \simeq_\ell Y$, if there exists some graph $Z \subseteq X, Y$ such that $\mathbb{P}_\ell(X) \cong \mathbb{P}_\ell(Y) \cong \mathbb{P}_\ell(Z)$. An $\ell$-path finite graph $G$ is said to be \emph{$\ell$-path minimal} if either $G$ is null, or for each $X \subset G$, $\mathbb{P}_\ell(X) \subset \mathbb{P}_\ell(G)$. Similar with Lemma \ref{lem_miniUniq}, we have:
\begin{propo}
$\simeq_\ell$ defines an equivalence relation on $\ell$-path finite graphs. Further, each $\ell$-path equivalence class contains a unique (up to isomorphism) minimal graph, which is $\ell$-path minimal. Moreover, an $\ell$-path minimal graph is a subgraph of every graph in its $\ell$-path equivalence class.
\end{propo}
%\pf
%Let $X \simeq_\ell Y$ and $Y \simeq_\ell Z$. By definition there are minimal subgraphs $Y_1$ and $Y_2$ of $Y$ such that $\mathbb{P}_\ell(X) \cong \mathbb{P}_\ell(Y_1)  = \mathbb{P}_\ell(Y) = \mathbb{P}_\ell(Y_2) \cong \mathbb{P}_\ell(Z)$. By the minimality, a unit of $Y$ belongs to $Y_i$ if and only if it belongs to an $\ell$-path of $Y$. So $Y_1 = Y_2$ is $\ell$-path minimal, and $X \simeq_\ell Z$. Therefore, $\simeq_\ell$ defines an equivalence relation, $Y_1$ is the unique $\ell$-path minimal graph in its class, and is a subgraph of every graph $X$ in the class.
%\qed

In the following, we exemplify that an $\ell$-minimal graph may not be $\ell$-path minimal. And an $\ell$-path minimal graph may not be an induced subgraph of another graph in its $\ell$-path equivalence class.
\begin{exa}
Let $G$ be a graph obtained from a path $P = [v_0,\ldots, v_4]$ by adding an edge $e$ between $v_1$ and $v_3$. By Lemmas \ref{lem_cycleIncident} and \ref{lem_minInci}, $G$ and $P$ are $4$-minimal. Since $\mathbb{P}_4(G) = \mathbb{P}_4(P) \cong K_1$, $G$ is not $\ell$-path minimal. Moreover, $P$ is $4$-path minimal and is a subgraph but not an induced subgraph of $G$.
\end{exa}

The example below indicates that, a graph can be both $\ell$-minimal and $\ell$-path minimal. But even so its $\ell$-link and $\ell$-path graphs may not be isomorphic.
\begin{exa}
The complete graph $K_{\ell + 1}$ is both $\ell$-minimal and $\ell$-path minimal. Clearly, for $\ell = 0, 1$ and $2$, $\mathbb{L}_\ell(K_{\ell + 1}) = \mathbb{P}_{\ell}(K_{\ell + 1})$ is isomorphic to $K_1, K_1$ and $K_3$ respectively. Now let $\ell \geqslant 3$. By \cite[Theorem 1.2]{JiaWood2013}, $\mathbb{L}_\ell(K_{\ell + 1})$ contains a $K_{\ell + 1}$-minor. However, $\mathbb{P}_{\ell}(K_{\ell + 1})$ consists of $\ell!/2$ disjoint cycles of length $\ell + 1$.
\end{exa}

%
%
%Let $s, \ell \geqslant 0$ be integers, $G$ be a graph.  For each $\ell$-path $L$ of $G$, \emph{$\mathcal {J}_G(L, s)$} be the set of
%pairs $(L, R)$ such that $R \in \mathscr{P}_s(G)$ is incident to $L$, and $j_G(L, s) := |\mathcal {J}_G(L, s)|$. Let \emph{$\mathcal {J}_G(\ell, s)$} be the union of $\mathcal {J}_G(L, s)$ over all $L \in \mathscr{P}_\ell(G)$ and \emph{$j_G(\ell, s):= |\mathcal {J}_G(\ell, s)|$}. Then $j_G(\ell, s) = \sum_{L \in \mathscr {P}_\ell(G)}j_G(L, s)$ and:
%\begin{lem}
%Let $G$ be a graph, $\ell \geqslant s \geqslant 0$, and $L \in \mathscr{P}_\ell(G)$. Then $j_G(L, s) = \ell - s + 1$, and $j_G(\ell, s) = (\ell - s + 1)|\mathscr{P}_\ell(G)|$.
%\end{lem}
%
%$R \in \mathscr{P}_s(G)$ is said to be \emph{$\ell$-path incident} if it is incident to some $L \in \mathscr{P}_\ell(G)$.
%
%\begin{lem}
%Let $\ell \geqslant 0$, $G$ be a graph. Then $G$ is $\ell$-path minimal iff each edge of $G$ is $\ell$-path incident. And in this case, each vertex of $G$ is also $\ell$-path incident. As a consequence, $j_G(\ell, 1) \geqslant \max\{\nu(G), \epsilon(G)\}$.
%\end{lem}

By definition $\mathbb{Q}_0(H) = \{H\}$ if $H$ is simple. Moreover $\mathbb{Q}_\ell(K_0) = \{K_0\}$, and $\mathbb{Q}_\ell(K_1)$ consists of an $\ell$-path. $\mathbb{Q}_1(K_2)$ consists of a $2$-path and a $2$-cycle. When $\ell \geqslant 2$, $\mathbb{Q}_\ell(K_2)$ consists of an $(\ell + 1)$-path. Similar with Lemma \ref{lem_lminiFinite}, the order, size and the total number of minimal $\ell$-path roots of a finite graph are bounded.
\begin{lem}\label{lem_pathrootosnbound}
Let $\ell \geqslant 1$ be an integer, and $H$ be a finite graph such that $n := n(H) \geqslant 2$. Let $c := c(H)$, and $G \in \mathbb{Q}_{\ell}(H)$. Then $n(G) \leqslant \ell n + c$ and $m(G) \leqslant \ell n$, with each equality holds if and only if  $G$ is a disjoint union of $\ell$-paths. Moreover, $\mathbb{Q}_{\ell}(H)$ contains at most $(\ell n + c)^{2\ell n}$ graphs, in which at most $\frac{3}{2}(\ell n + 1)^{\ell n - 1}$ are trees, and at most $\frac{3}{4}c(c + 1)(\ell n + c)^{\ell n - 1}$ are forests.
\end{lem}

\bibliographystyle{plain}
\bibliography{Xbib}
%\tableofcontents

\end{document}